\documentclass[12pt]{article}
\usepackage{amsmath,amscd,amsthm,amssymb}

\usepackage[russian]{babel}
\begin{document}

\begin{center}
\LARGE Some aspects of harmonic analysis related to Gegenbauer expansions on the half-line
\end{center}

\

\centerline{\large Vagif S. Guliyev and Elman J. Ibrahimov}

\

\

\textbf{Abstract} In this paper we consider the generalized shift operator, generated by the Gegenbauer differential operator 
\[
G =\left(x^2-1\right)^{\frac{1}{2}-\lambda } \frac{d}{dx} \left(x^2-1\right)^{\lambda+\frac{1}{2}}\frac{d}{dx}.
\]
Maximal function ($ G- $ maximal function), generated by the Gegenbauer differential operator $ G $ is investigated. 
The $ L_{p,\lambda} $ -boundedness for the $ G- $ maximal function is obtained. The concept of potential of Riesz-Gegenbauer 
is introduced and for it the theorem of Sobolev type is proved.

\


\textbf{Mathematics Subject Classifications (2000) Primary } 42B20, 42B25, 42B35

\

\

\textbf{Introduction} The Hardy-Lettlewood maximal function is an important tool of harmonic analysis. It was first introduced by Hardy and Littlewood in 1930 (see [21]) for $ 2\pi- $ periodical functions, and later it was extendet to the Euclidean spaces, some weighted measure spaces (see [4, 28, 29]) , symmetric spaces (see [5, 26]), various Lie groups [9], for the Jacobi-type hypergroups [6, 7], for Chebli-Trimeche hypergroups [1], for the one-dimensional Bessel-Kingman hypergroups [27], for the $n-$ dimensional Bessel-Kingman hypergroups $\left( n\geq1\right) $ [10, 11, 13], for Morrey-Bessel spaces [2, 3, 12, 14], for Laguerre hypergroup [15, 16, 22, 25].
The structure of the paper is as follows. In Section 1 we present some definitions, notation and auxiliary results. In Section 2 the $ L_{p,\lambda} $ boundedness of the $ G -$ maximal function is proved. In Section 3 we introduce and study some embeddings into the function spaces, associated with the Gegenbauer differential operator. In Section 4 we introduce definition Riesz-Gegenbaur potential and is for to Sobolev type theorem is proved.

\section{DEFINITIONS, NOTATION AND AUXILIARY RESULTS}

Let $H(x,r)=(x-r, x+r)\cap [0, \infty), r \in (0, \infty), x\in [0, \infty)$. For all measurable set
$E\subset [0, \infty)~$    $\mu E\equiv |E|_{\lambda}= \int \limits_{E} sh^{2\lambda }\,\,t dt$. For $1\leq p\leq \infty$ let $L_{p}([0, \infty), G)\equiv L_{p,\lambda}[0, \infty)$ be the space of functions measurable on $[0, \infty)$ with the finite norm

\[\Vert f \Vert_{L_{p,\lambda}}=\left( \int_{0}^{\infty}|f(ch\,\,t)|^{p} sh^{2\lambda }\,\,t dt\right)^{\frac{1}{p}}, 1\leq p< \infty,\]

\[\Vert f \Vert_{\infty, \lambda}= \operatorname*{ess\,sup}\limits_{t\in [0, \infty )} |f(ch\,\,t)|,\,\,\, p=\infty\]

 Analogy by [7] we define Gegenbauer maximal functions is as follows:

\[
M_{G} f\left(ch\,\,x\right)=\mathop{\sup }\limits_{r>0} \frac{1}{\left|H(0,r\right|_{\lambda } } \int\limits _{0}^{r}A_{ch\,\,t}^{\lambda } \left|f\left(ch\,\,x\right)\right| d\mu \left(t\right),
\]
\[
M_{\mu } f\left(ch\,\,x\right)=\mathop{\sup }\limits_{r>0} \frac{1}{\left|H\left(x,r\right)\right|} _{\lambda } \int\limits _{H(x,r)}\left|f(ch\,\,t\right| d\mu \left(t\right),\,\, d\mu(t)=sh^{2\lambda }\,\,tdt,
\]
\[
\left|H(0,r)\right|_{\lambda } =\int\limits _{0}^{r}sh^{2\lambda }\,\,tdt, \,\,\,\,  \left|H(x,r)\right|_{\lambda } =\int\limits _{H(x,r)}sh^{2\lambda } \,\,tdt.
\]

Here
\[
A_{ch\,\,t}^{\lambda } f(ch\,\,x)=\frac{\Gamma (\lambda +\frac{1}{2} )}{\Gamma (\frac{1}{2} )\Gamma \left(\lambda \right)} \int\limits _{0}^{\pi }f(ch\,\,xch\,\,t-sh\,\,xsh\,\,t\cos\,\, \varphi  )(\sin\,\,\varphi )^{2\lambda -1} d\varphi
\]
denote the generalized shift operator, associated with the Gegenbauer differential operator
\[
G =\left(x^2-1\right)^{1/2-\lambda } \frac{d}{dx} \left(x^2-1\right)^{\lambda+1/2}\frac{d}{dx}.
\]

By $A \lesssim B$ we mean that $A \le C B$ with some positive constant $C$ independent of appropriate quantities. If $A \lesssim B$ and $B \lesssim A$, we
write $A\approx B$ and say that $A$ and $B$ are  equivalent.

Further we'll need some auxiliary assertions.

\textbf{Lemma 1}  {\it For   $0<\lambda <{1 \mathord{\left/{\vphantom{1 2}}\right.\kern-\nulldelimiterspace} 2} $ the following correlations are true:
\[
\left|H(0,r)\right|_{\lambda } \sim \left\{\begin{array}{l} {(sh\,\,\frac{r}{2} )^{2\lambda +1} ,0<r\le c;\, \, \ \ \ \, \, \, \, \, \, \, \, \, \, \, \, \, (a)} \\ \\{(ch\,\,\frac{r}{2} )^{4\lambda },\quad c<r<\infty ,\, \, \, \, \, \, \, \, \, \, \, \, \, \, \, \, \, \, (b)}
\end{array}\right.
\]
where $c$ denotes positive constant.}


 \textbf{Proof}. Let first $0<r\le c,$ then
\[
\left|H(0,r)\right|_{\lambda } =\int\limits _{0}^{r}sh^{2\lambda } \,\,tdt=\int\limits _{0}^{r}(sh\,\,t)^{2\lambda -1}  d(ch\,\,t)
=\int\limits _{0}^{r}(ch^{2}\,\, t-1)^{\lambda -\frac{1}{2} }  d(ch\,\,t)
\]
\[=\int\limits _{1}^{ch\,\,r}(t-1)^{\lambda -\frac{1}{2} }  (t+1)^{\lambda -\frac{1}{2} } dt
\ge (ch\,\,r+1)^{\lambda -\frac{1}{2} } \int\limits _{1}^{ch\,\,r}(t-1)^{\lambda -\frac{1}{2} }  dt
\]
\[
\ge (ch1+1)^{\lambda -\frac{1}{2} } \frac{(t-1)^{\lambda +\frac{1}{2} } }{\lambda +\frac{1}{2} } \left|_{1}^{ch\,\,r} \right. =\frac{2(ch\,\,r-1)^{\lambda +\frac{1}{2} } }{(2\lambda +1)(1+ch1)^{\frac{1}{2} -\lambda } }
\]
\[
 =\frac{2^{2\lambda +2} }{(2\lambda +1)(1+ch1)^{\frac{1}{2} -\lambda } } (sh\,\,\frac{r}{2} )^{2\lambda +1} .
\eqno(1)
\]

 On the other hand,
\[
\left|H(0,r)\right|_{\lambda } =\int\limits _{0}^{r}sh^{2\lambda } \,\,tdt=\int\limits _{1}^{ch\,\,r}(t-1)^{\lambda -\frac{1}{2} }  (t+1)^{\lambda -\frac{1}{2} } dt\le 2^{\lambda -\frac{1}{2} } \int\limits _{1}^{ch\,\,r}(t-1)^{\lambda -\frac{1}{2} }  dt
\]
\[
=\frac{2^{\lambda +\frac{1}{2} } }{2\lambda +1} (t-1)^{\lambda +\frac{1}{2} } \left|_{1}^{ch\,\,r} \right. =\frac{2^{\lambda +\frac{1}{2} } }{2\lambda +1} (ch\,\,r-1)^{\lambda +\frac{1}{2} } =\frac{2^{2\lambda +1} }{2\lambda +1} \left(sh\,\,\frac{r}{2} \right)^{2\lambda +1} .
\eqno(2)
\]

Let, now $c<r<\infty .$ Then
$$
\left|H(0,r)\right|_{\lambda } =\int\limits _{0}^{r}sh^{2\lambda } \,\,tdt=\int\limits _{0}^{r}(sh\,\,t)^{2\lambda -1}  d(ch\,\,t)=\int\limits _{0}^{r}(ch^{2}\,\,  t-1)^{\lambda -\frac{1}{2} } d(ch\,\,t)
$$
$$
=\int\limits _{1}^{ch\,\,r}\frac{(t-1)^{\lambda -\frac{1}{2} } }{(t+1)^{\frac{1}{2} -\lambda } }  dt
 \ge (ch\,\,r+1)^{\lambda -\frac{1}{2} } \int\limits _{1}^{ch\,\,r}(t-1)^{\lambda -\frac{1}{2} }  dt
 $$
 $$
 =\left.(ch\,\,r+1)^{\lambda -\frac{1}{2} } \frac{(t-1)^{\lambda +\frac{1}{2} } }{\lambda +\frac{1}{2} } \right|_{1}^{ch\,\,r}  =\frac{2}{2\lambda +1} \frac{(ch\,\,r-1)^{\lambda +\frac{1}{2} } }{(ch\,\,r+1)^{\frac{1}{2} -\lambda } }
$$
$$
=\frac{2}{2\lambda +1} \frac{\left(2sh^{2}\,\, \frac{r}{2} \right)^{\lambda +\frac{1}{2} } }{\left(2ch^{2}\,\, \frac{r}{2} \right)^{\frac{1}{2} -\lambda } } =\frac{2^{2\lambda +2} }{(2\lambda +1)2^{1-2\lambda } } \frac{\left(sh\,\,\frac{r}{2} \right)^{2\lambda +1} }{\left(ch\,\,\frac{r}{2} \right)^{1-2\lambda } }
$$
$$
\ge \frac{2^{4\lambda +1} }{(2\lambda +1)3^{2\lambda +1} } \left(ch\,\,\frac{r}{2} \right)^{4\lambda }
\Leftrightarrow \frac{4^{2\lambda +1} }{2\lambda +1} \left(3sh\,\,\frac{r}{2} \right)^{2\lambda +1}
$$
$$
\ge \frac{4^{2\lambda +1} }{2\lambda +1} \left(ch\,\,\frac{r}{2} \right)^{2\lambda +1} \Leftrightarrow 3sh\,\,\frac{r}{2} \ge ch\,\,\frac{r}{2} \Leftrightarrow 3\frac{e^{{r \mathord{\left/{\vphantom{r 2}}\right.\kern-\nulldelimiterspace} 2} } -e^{-{r \mathord{\left/{\vphantom{r 2}}\right.\kern-\nulldelimiterspace} 2} } }{2}
$$
$$
 {\ge \frac{e^{{r \mathord{\left/{\vphantom{r 2}}\right.\kern-\nulldelimiterspace} 2} } +e^{-{r \mathord{\left/{\vphantom{r 2}}\right.\kern-\nulldelimiterspace} 2} } }{2} \Leftrightarrow 3(e^{r} -1)\ge e^{r} +1\Leftrightarrow 2e^{r} \ge 4,}
$$
  which takes place for $r\ge c\ge 1$.

So,
$$
\left|H(0,r)\right|_{\lambda } \ge \frac{2^{4\lambda +1} }{(2\lambda +1)3^{2\lambda +1} } \left(ch\,\,\frac{r}{2} \right)^{4\lambda } .
\eqno(3)
$$

Estimate above $\left|H(0,r)\right|_{\lambda } $.
\[
\left|H(0,r)\right|_{\lambda } =\int\limits _{0}^{r}sh^{2\lambda } \,\,tdt=\int\limits _{0}^{r}\left(2sh\,\,\frac{t}{2}  ch\,\,\frac{t}{2} \right)^{2\lambda } dt
\]
\[
=2^{2\lambda +1} \int\limits _{0}^{r}\left(sh\,\,\frac{t}{2}  \right)^{2\lambda } \left(ch\,\,\frac{t}{2} \right)^{2\lambda -1} d\left(sh\,\,\frac{t}{2} \right)\le 2^{2\lambda +1} \int\limits _{0}^{r}\left(sh\,\,\frac{t}{2}  \right)^{4\lambda -1} d\left(sh\,\,\frac{t}{2} \right)
\]
\[\left.=\frac{2^{2\lambda +1} }{4\lambda } \left(sh\,\,\frac{t}{2} \right)^{4\lambda } \right|_{0}^{r}  =\frac{4^{\lambda } }{2\lambda } \left(sh\,\,\frac{r}{2} \right)^{4\lambda } \le \frac{4^{\lambda } }{2\lambda } \left(ch\,\,\frac{r}{2} \right)^{4\lambda } .\eqno\left(4\right)
\]

Combine (1)-(4), we obtain assertion of lemma 1.

\textbf{Lemma 2} {\it Let $0<\lambda <1/2$ and $x\in [0,\infty)$, $r\in (0,\infty)$. Then the following estimates are true: for $0<r\le c$
$$
\left|H(x,r)\right|_{\lambda } \le c_{\lambda } \left\{\begin{array}{l} {r^{2\lambda +1} ,\,\,\, 0\le x\le r\le c;} \\ \ \\ {ch^{2\lambda }\,\, x,\, \, \, \, \, r<x<\infty\ (r\le c<x<\infty) .} \end{array}\right.                                  \eqno(a)
$$

And for $c<r<\infty $}
$$
\left|H(x,r)\right|_{\lambda } \le c_{\lambda }
\left\{\begin{array}{l} {ch^{2\lambda }\,\, r,\,\,\, 0<x\le 2r \  (0<x\le 2c<2r);} \\ \ \\ {ch^{2\lambda }\,\, xch^{2\lambda }\,\, r,\,\,\,2r<x<\infty \ (2c<2r<x<\infty).}
\end{array}\right.                             \eqno(b)
$$

Here and further $c_{\lambda } ,c_{\alpha ,\lambda } ,c_{\alpha ,\lambda ,p} $ will denote some constants, depending only on subscribed indexes and generally speaking different  in different formulas.

\textbf{Proof.} First we consider case when $0<r\le c$ and $x\in [0,\infty).$

Let $0\le t\le 2c,$ Then we have
$$
t\le sh\,\,t\le e^{2c} t.\eqno(5)
$$

We prove left-hand part of this estimate. We consider the function\linebreak $f(t)=sh\,\,t-t$. As, $f'(t)=ch\,\,t-t\ge 0$, then $f(t)$ increases in $[0,\infty)$, and that takes the smallest valuer for $t=0$, $f(0)=0$, consequently  $f(t)\ge 0$ equivalent to $sh\,\,t\ge t$.

We prove right-hand part of estimate (5).
\[
\frac{e^{t} -e^{-t} }{2} \le e^{2c} \cdot t\, \, \, \, \Leftrightarrow e^{2t} -1\le 2\cdot e^{2c+t} \cdot t\, \, \Leftrightarrow e^{2t} \le 2\cdot e^{2c+t} \cdot t+1.
\]

We consider the function $f(t)=2\cdot e^{2c+t} \cdot t+1-e^{2t} $.
\[
f'(t)=2\cdot e^{2c+t} +2\cdot e^{2c+t} \cdot t-2e^{2t} =2e^{t} \left(e^{2c} +t\cdot e^{2c} -e^{t} \right)
\]
$$
\ge e^{2c} \left(t+1\right)-e^{t}\ge e^{2c} -e^{t} \ge 0,\ \mbox{as,}\ \ t\le 2c.
$$

Thus, the estimate (5) is proved.

Hence follows that for $0\le x\le r\le c$
$$
\left|H(x,r)\right|_{\lambda } =\int\limits _{0}^{x+r}sh^{2\lambda } \,\,tdt\le e^{2c} \int\limits _{0}^{2r}t^{2\lambda }  dt=\frac{e^{2c} \cdot 2^{2\lambda } }{2\lambda +1} \cdot r^{2\lambda +1} . \eqno(6)
$$%

 $r<x\le c$
$$
\left|H(x,r)\right|_{\lambda } =\int\limits _{x-r}^{x+r}sh^{2\lambda } \,\,tdt\le e^{2c} \int\limits _{x-r}^{x+r}t^{2\lambda }  dt\le e^{2c} \cdot r\cdot (x+r)^{2\lambda }
$$$$
\le e^{2c} \cdot r\cdot (2x)^{2\lambda } \le c_{\lambda } ch^{2\lambda }\,\, x.
\eqno(7)
$$%

Let, now $0<r\le c\le x<\infty $, then we have
$$
\left|H(x,r)\right|_{\lambda } =\int\limits _{x-r}^{x+r}sh^{2\lambda } \,\,tdt\le 2r\cdot sh^{2\lambda } (x+r)=2r(sh\,\,xch\,\,r+ch\,\,xsh\,\,r)^{2\lambda }
$$
$$
\le 2r(sh\,\,xch\,\,c+ch\,\,xsh\,\,c)^{2\lambda }
 \le 2r(2ch\,\,xch\,\,c)^{2\lambda } \le c_{\lambda } ch^{2\lambda }\,\, x.
 \eqno\left(8\right)
$$

Now, we consider case, when  $c<r<\infty ,\,\,\,x\in \left[0,\infty \right).$

Let, $0<x\le 2c<2r$. Enterling as with proof of the estimate (4), we obtain
\[
\left|H(x,r)\right|_{\lambda } =\int\limits _{0}^{x+r}sh^{2\lambda }\,\,tdt=\frac{4^{\lambda } }{2\lambda }  sh^{4\lambda } \frac{t}{2} \left|_{0}^{x+r} \right. =\frac{4^{\lambda } }{2\lambda } sh^{4\lambda } \,\, \,\, \frac{x+r}{2}
 \]
 \[
 =\frac{4^{\lambda } }{2\lambda } \left(sh\,\,\frac{x}{2} ch\,\,\frac{r}{2} +ch\frac{x}{2} sh\,\,\frac{r}{2} \right)^{4\lambda }   \le c_{\lambda } (sh\,\,c\, ch\,\,\frac{r}{2} +ch\,\,c\, sh\,\,\frac{r}{2} )^{4\lambda } \le c_{\lambda } ch^{4\lambda } \,\, \frac{r}{2} .
\eqno(9)
\]

Let, now $2c<2r<x<\infty ,$ then
\[
 \left|H(x,r)\right|_{\lambda } =\int\limits _{x-r}^{x+r}sh^{2\lambda } \,\,tdt\le \left.\frac{4^{\lambda } }{2\lambda }  \,\, \,\,\frac{t}{2} \right|_{x-r}^{x+r}
  =\frac{4^{\lambda } }{2\lambda } \left(sh^{4\lambda } \,\, \frac{x+r}{2} -sh^{4\lambda } \,\, \frac{x-r}{2}\right)
 \]
 \[
  \le \frac{4^{\lambda } }{2\lambda } sh^{4\lambda } \,\, \frac{x+r}{2}   \le c_{\lambda } ch^{4\lambda } \,\, \frac{x}{2} ch^{4\lambda } \,\, \frac{r}{2} \le c_{\lambda } ch^{2\lambda }\,\, xch^{2\lambda }\,\, r.
\eqno\left(10\right)
\]

From (9)  and  (10) follows that at  $c<r<\infty $ and  $0<x<\infty $
\[
\left|H(x,r)\right|_{\lambda } \le c_{\lambda } \left\{\begin{array}{l} {ch^{2\lambda }\,\, r,\, \, \, \, \, \, \, \, \, \, \, \, \, \, \, \, \, \, \, \, \, 0<x\le 2c<2r;\qquad\qquad \left(11\right)}\\ \ \\ {ch^{2\lambda }\,\, xch^{2\lambda }\,\, r,\ \  2c<2r<x<\infty .\qquad\qquad\left(12\right)}
\end{array}\right.
\]

Assertion of lemma 2 follows from (6)-(8), (11) and (12).

\section{$L_{p,\lambda}$ -boundedness of $G$- maximal operator}

\textbf{Theorem 1} {\it For $0\le x<\infty $ and $0<r<\infty $ the inequality is valid
\[
M_{G} f(ch\,\,x)\le c_{\lambda } M_{\mu } f(ch\,\,x),
\]
where $c_{\lambda } $ is some positive constant.}

\textbf{Proof.} Consider the integral
\[
{I(x,r)=\int\limits _{0}^{r}A_{ch\,\,t}^{\lambda }  \left|f(ch\,\,x)\right|sh^{2\lambda }\,\,tdt}
\]
\[ {=\frac{\Gamma (\lambda +\frac{1}{2} )}{\Gamma (\frac{1}{2} )\Gamma (\lambda )} \int\limits _{0}^{r}\left\{\int\limits _{0}^{\pi }\left|f(ch\,\,x\cdot ch\,\,t-sh\,\,x\cdot sh\,\,t\cos\,\, \varphi  \right|(\sin\,\,\varphi )^{2\lambda -1} d\varphi  \right\} sh^{2\lambda }\,\,tdt.}
\]

Making in internal integral replacing\\ $z=ch\,\,x\cdot ch\,\,t-sh\,\,x\cdot sh\,\,t\cos\,\, \varphi  ,$ we get that
\[
{\cos\,\, \varphi  =\frac{ch\,\,x\cdot ch\,\,t-z}{sh\,\,x\cdot sh\,\,t} ,\varphi =\arccos \frac{ch\,\,x\cdot ch\,\,t-z}{sh\,\,x\cdot sh\,\,t} ,}
\]
\[
d\varphi =\frac{dz}{\sqrt{1-(\frac{ch\,\,x\cdot ch\,\,t-z}{sh\,\,x\cdot sh\,\,t} )^{2} } sh\,\,x\cdot sh\,\,t}
\]
\[
=(sh^{2}\,\, x\cdot sh^{2}\,\, t-ch^{2}\,\, x\cdot ch^{2}\,\, t+2\cdot z\cdot ch\,\,x\cdot ch\,\,t-z^{2} )^{-\frac{1}{2} } dz.
\]

As,
\[
sh^{2}\,\, x\cdot sh^{2}\,\, t-ch^{2}\,\, x\cdot ch^{2}\,\, t
\]\[
=(ch^{2}\,\, x-1)sh^{2}\,\, t-ch^{2}\,\, x\cdot ch^{2}\,\, t=ch^{2}\,\, x\cdot sh^{2}\,\, t-sh^{2}\,\, t-ch^{2}\,\, x\cdot ch^{2}\,\, t
\]
\[ =-sh^{2}\,\, t+ch^{2}\,\, x(sh^{2}\,\, t-ch^{2}\,\, t)=-sh^{2}\,\, t-ch^{2}\,\, x,
\]
that
\[
d\varphi =(2z\cdot ch\,\,x\cdot ch\,\,t-sh^{2}\,\, t-ch^{2}\,\, x-z^{2} )^{-\frac{1}{2} } dz
\]
and
\[
(\sin\,\,\varphi )^{2\lambda -1} =(2z\cdot ch\,\,x\cdot ch\,\,t-sh^{2}\,\, t-ch^{2}\,\, x-z^{2} )^{\lambda -\frac{1}{2} } (sh\,\,x\cdot sh\,\,t)^{1-2\lambda } .
\]

Then $I(x,r)$ makes a list of form
$$
I(x,r)=\frac{\Gamma (\lambda +\frac{1}{2} )}{\Gamma (\frac{1}{2} )\Gamma (\lambda )}
$$$$
\times \int\limits _{0}^{r}\left\{\int\limits _{ch\,\,(x-t)}^{ch\,\,(x+t)}\left|f(z)\right|(2z\cdot ch\,\,x\cdot ch\,\,t-sh^{2}\,\, t-ch^{2}\,\, x-z^{2} )^{\lambda -1} (sh\,\,x)^{1-2\lambda } dz \right\} sh\,\,tdt.
\eqno(13)
$$

Transform expansion
\[
2z\cdot ch\,\,x\cdot ch\,\,t-sh^{2}\,\, t-ch^{2}\,\, x-z^{2}
\]\[
 =2z\cdot ch\,\,x\cdot ch\,\,t-sh^{2}\,\, t(ch^{2}\,\, x-sh^{2}\,\, x)-ch^{2}\,\, x-z^{2}
\]
\[ =2z\cdot ch\,\,x\cdot ch\,\,t-sh^{2}\,\, t\cdot ch^{2}\,\, t+sh^{2}\,\, t\cdot sh^{2}\,\, x-ch^{2}\,\, x-z^{2}
\]
\[=2z\cdot ch\,\,x\cdot ch\,\,t+sh^{2}\,\, t\cdot sh^{2}\,\, x-(ch^{2}\,\, t-1)ch^{2}\,\, x-ch^{2}\,\, x-z^{2}
\]
\[
=2z\cdot ch\,\,x\cdot ch\,\,t+sh^{2}\,\, x\cdot (ch^{2}\,\, t-1)-ch^{2}\,\, t\cdot ch^{2}\,\, x-z^{2} (ch^{2}\,\, x-sh^{2}\,\, x)
\]
\[
=2z\cdot ch\,\,x\cdot ch\,\,t+sh^{2}\,\, x\cdot ch^{2}\,\, t-sh^{2}\,\, x-ch^{2}\,\, t\cdot ch^{2}\,\, x
-z^{2} \cdot ch^{2}\,\, x-z^{2} \cdot sh^{2}\,\, x
\]\[
=2z\cdot ch\,\,x\cdot ch\,\,t-sh^{2}\,\, x-ch^{2}\,\, t-z^{2} ch^{2}\,\, x-z^{2} sh^{2}\,\, x
=(z^{2} -1)sh^{2}\,\, x-(ch\,\,t-z\cdot ch\,\,x)^{2}
 \]
\[
=\left(z^{2} -1\right)sh^{2}\,\, x\left[1-\left(\frac{ch\,\,t-z\cdot ch\,\,x}{\sqrt{z^{2} -1} \cdot sh\,\,x} \right)^{2} \right].
\eqno \left(14\right)
\]

Taking into account (13) and (14) we get
$$
I(x,r)=\frac{\Gamma (\lambda +\frac{1}{2} )}{\Gamma (\frac{1}{2} )\Gamma (\lambda )} \int\limits _{0}^{r}\left\{\int\limits _{ch\,\,(x-t)}^{ch\,\,(x+t)}\left|f(z)\right|(z^{2} -1)^{\lambda -1}\right.
$$
$$\times\left. \left[1-\left(\frac{ch\,\,t-z\cdot ch\,\,x}{\sqrt{z^{2} -1} \cdot sh\,\,x} \right)^{2} \right] ^{\lambda -1} dz\right\}\frac{sh\,\,t}{sh\,\,x} dt.
\eqno(15)
$$%

Note that
\[
\frac{sh\,\,t}{sh\,\,x} =(z^{2} -1)^{\frac{1}{2} } \frac{\partial }{\partial t} \left(\frac{ch\,\,t-z\cdot ch\,\,x}{\sqrt{z^{2} -1} \cdot sh\,\,x} \right),
\]
rewrite (15) of form
$$
I(x,r)=\frac{\Gamma (\lambda +\frac{1}{2} )}{\Gamma (\frac{1}{2} )\Gamma (\lambda )} \int\limits _{0}^{r}\left\{\int\limits _{ch\,\,(x-t)}^{ch\,\,(x+t)}\left|f(z)\right|(z^{2} -1)^{\lambda -\frac{1}{2} } \right.
 $$
 $$
\left.\times\left[1-\left(\frac{ch\,\,t-z\cdot ch\,\,x}{\sqrt{z^{2} -1} \cdot sh\,\,x} \right)^{2} \right]^{\lambda -1}  \frac{\partial }{\partial t} \left(\frac{ch\,\,t-z\cdot ch\,\,x}{\sqrt{z^{2} -1} \cdot sh\,\,x} \right)\right\}dzdt.\eqno(16)
$$%

As, $ch\,\,(x-t)\le z\le ch\,\,(x+t),then$
$$
\left\{\begin{array}{l} {ch\,\,(x-r)\le z\le ch\,\,x} \\ {x-arcchz\le t\le r} \end{array}\right. \ \     \mbox{and}\    \ \left\{\begin{array}{l} {ch\,\,x\le z\le ch\,\,(x+r)} \\ {arcchz-x\le t\le r} \end{array}\right.
$$

And that's why changing the order of integration in (16) , we get
$$
I(x,r)=\frac{\Gamma (\lambda +\frac{1}{2} )}{\Gamma (\lambda )\Gamma (\frac{1}{2} )} \left(\int\limits _{ch\,\,(x-r)}^{ch\,\,x}dz \int\limits _{x-arcchz}^{r}dt +\int\limits _{ch\,\,x}^{ch\,\,(x+r)}dz \int\limits _{arcchz-x}^{r}dt \right). \eqno(17)
$$%

Consider the integral

\[
A(x,z,r)\equiv A(x,r)
\]
\[
=\int\limits _{x-arcchz}^{r}\left[1-\left(\frac{ch\,\,t-z\cdot ch\,\,x}{\sqrt{z^{2} -1} \cdot sh\,\,x} \right)^{2} \right]^{\lambda -1} \frac{\partial }{\partial t} \left(\frac{ch\,\,t-z\cdot ch\,\,x}{\sqrt{z^{2} -1} \cdot sh\,\,x} \right)dt .
\]

Putting here $u=\dfrac{ch\,\,t-z\cdot ch\,\,x}{\sqrt{z^{2} -1} \cdot sh\,\,x} ,$ we get
$$
A(x,z,r)\equiv A(x,r)=\int\limits _{-1}^{\frac{ch\,\,r-z\cdot ch\,\,x}{\sqrt{z^{2} -1} \cdot sh\,\,x} }(1-u^{2} )^{\lambda -1} du  \eqno(18)
$$%

On the power of even  $ch\,\,t$
$$
B(x,r)=\int\limits _{arcchz-x}^{r}\left[1-\left(\frac{ch\,\,t-z\cdot ch\,\,x}{\sqrt{z^{2} -1} \cdot sh\,\,x} \right)^{2} \right]^{\lambda -1} \frac{\partial }{\partial t} \left(\frac{ch\,\,t-z\cdot ch\,\,x}{\sqrt{z^{2} -1} \cdot sh\,\,x} \right)dt
$$
$$
=\int\limits _{-1}^{\frac{ch\,\,r-z\cdot ch\,\,x}{\sqrt{z^{2} -1} \cdot sh\,\,x} }(1-u^{2} )^{\lambda -1} du.\eqno(19)
$$%

Taking into account (18) and (19) in (17), we have
$$
I(x,r)=\frac{\Gamma (\lambda +\frac{1}{2} )}{\Gamma (\frac{1}{2} )\Gamma (\lambda )} \int\limits _{ch\,\,(x-r)}^{ch\,\,(x+r)}\left|f(z)\right|(z^{2} -1)^{\lambda -\frac{1}{2} } \int\limits _{-1}^{\frac{ch\,\,r-z\cdot ch\,\,x}{\sqrt{z^{2} -1} \cdot sh\,\,x} }(1-u^{2}  )^{\lambda -1} dudz . \eqno(20)
$$%

As, $ch\,\,(x-r)\le z\le ch\,\,(x+r),then$
$$
 \frac{ch\,\,r-z\cdot ch\,\,x}{\sqrt{z^{2} -1} \cdot sh\,\,x} \ge \frac{ch\,\,r-ch\,\,x\cdot ch\,\,(x+r)}{sh\,\,x\cdot sh(x+r)}
  =\frac{2ch\,\,r-2ch\,\,x\cdot ch\,\,(x+r)}{2sh\,\,x\cdot sh(x+r)}
  $$$$
  =\frac{2ch\,\,r-ch(2x+r)-ch\,\,r}{ch(2x+r)-ch\,\,r} =\frac{ch\,\,r-ch(2x+r)}{ch(2x+r)-ch\,\,r} =-1
 \eqno(21)
$$%

On the other hand for $ch\,\,(x-r)\le z\le ch\,\,(x+r),$
$$
\frac{ch\,\,r-z\cdot ch\,\,x}{\sqrt{z^{2} -1} \cdot sh\,\,x} \le \frac{ch\,\,r-ch\,\,x\cdot ch\,\,(x-r)}{sh\,\,x\left|sh(x-r)\right|} =\frac{2ch\,\,r-2ch\,\,x\cdot ch\,\,(x-r)}{2sh\,\,x\cdot sh(r-x)}
 $$$$
 =\frac{2ch\,\,r-ch(2x-r)-ch\,\,r}{ch\,\,r-ch(2x-r)} =\frac{ch\,\,r-ch(2x-r)}{ch\,\,r-ch(2x-r)} =1.
\eqno(22)
$$%

From (21) and (22) follows that for $ch\,\,(x-r)\le z\le ch\,\,(x+r),$ and $0<x\le r$
$$
-1\le \frac{ch\,\,r-z\cdot ch\,\,x}{\sqrt{z^{2} -1} \cdot sh\,\,x} \le 1.\eqno(23)
$$%

From (23) follows that for $0<x\le r\le c$
$$
A(x,r)=\int\limits _{-1}^{\frac{ch\,\,r-z\cdot ch\,\,x}{\sqrt{z^{2} -1} \cdot sh\,\,x} }(1-u^{2} )^{\lambda -1} du\le \int\limits _{-1}^{1}(1-u^{2} )^{\lambda -1} du=\frac{\Gamma (\frac{1}{2} )\Gamma (\lambda )}{\Gamma (\lambda +\frac{1}{2} )}   . \eqno(24)
$$%

But then taking into account (24) and (20), we obtain that for $0<x\le r\le c$
$$
I(x,r)\le \int\limits _{ch\,\,(x-r)}^{ch\,\,(x+r)}\left|f(z)\right|(z^{2} -1)^{\lambda -\frac{1}{2} }  dz=\int\limits _{x-r}^{x+r}\left|f(ch\,\,t)\right|sh^{2\lambda }\,\,tdt . \eqno(25) $$%

Now, let $c<r<x<\infty $ and $ch\,\,(x-r)\le z\le ch\,\,(x+r)$\\ $(c<r<x<\infty).$

Then we have
$$
\frac{ch\,\,r-z\cdot ch\,\,x}{\sqrt{z^{2} -1} \cdot sh\,\,x}\! \le\! \frac{ch\,\,x-z\cdot ch\,\,x}{\sqrt{z^{2} -1} sh\,\,x} \!=\!\frac{(1-z)ch\,\,x}{\sqrt{z^{2} -1}  sh\,\,x} \!=-\frac{\sqrt{z-1} ch\,\,x}{\sqrt{z+1} sh\,\,x} \le 0.
\eqno(26)
$$%

From (21) follows, that
\[
\max (1-u)^{\lambda -1} \le \mathop{\max }\limits_{-1\le u\le 0} (1-u)^{\lambda -1} =\max (2^{\lambda -1} ,1)=1,
\]
\[
-1\le u\le \frac{ch\,\,r-z\cdot ch\,\,x}{\sqrt{z^{2} -1} \cdot sh\,\,x} .
\]

Taking into account this circumstance, for the integral $A(x,r)$ we obtain of (18)
$$
A(x,r)=\int\limits _{-1}^{\frac{ch\,\,r-z\cdot ch\,\,x}{\sqrt{z^{2} -1\cdot sh\,\,x} } }\left(1-u^{2} \right)^{\lambda -1} du
$$$$
\le \int\limits _{-1}^{\frac{ch\,\,r-z\cdot ch\,\,x}{\sqrt{z^{2} -1\cdot sh\,\,x} } }\left(1+u\right)^{\lambda -1} du=\left.\frac{1}{\lambda } \left(1+u\right)^{\lambda } \right|_{-1} ^{\frac{ch\,\,r-z\cdot ch\,\,x}{\sqrt{z^{2} -1\cdot sh\,\,x} } }  =\frac{1}{\lambda } \left(1+\frac{ch\,\,r-z\cdot ch\,\,x}{\sqrt{z^{2} -1\cdot sh\,\,x} } \right)^{\lambda }
$$
$$
=\frac{1}{\lambda } \left(1-\frac{z\cdot ch\,\,x-ch\,\,r}{\sqrt{z^{2} -1} \cdot sh\,\,x} \right)^{\lambda } \le \frac{1}{\lambda } \left[1-\left(\frac{z\cdot ch\,\,x-ch\,\,r}{\sqrt{z^{2} -1} \cdot sh\,\,x} \right)^{2} \right]^{\lambda } .
\eqno(27)
$$%

We find extremum of the function
\[
f(z)=1-\left(\frac{z\cdot ch\,\,x-ch\,\,r}{\sqrt{z^{2} -1} \cdot sh\,\,x} \right)^{2} .
\]
\[
 f'(z)=-2\left(\frac{z\cdot ch\,\,x-ch\,\,r}{\sqrt{z^{2} -1} \cdot sh\,\,x} \right)
 \]
 \[
 \times\frac{(z^{2} -1)sh\,\,x\cdot ch\,\,x-z^{2} sh\,\,x\cdot ch\,\,x+z\cdot ch\,\,r\cdot sh\,\,x}{(z^{2} -1)^{\frac{3}{2} } sh^{2}\,\, x}
 \]
 \[ =-2\left(\frac{z\cdot ch\,\,x-ch\,\,r}{\sqrt{z^{2} -1} \cdot sh\,\,x} \right)\frac{z\cdot ch\,\,r\cdot sh\,\,x-ch\,\,x\cdot sh\,\,x}{(z^{2} -1)^{\frac{3}{2} } sh^{2}\,\, x}
 \]
 \[
  =\frac{2(z\cdot ch\,\,x-ch\,\,r)(ch\,\,x-z\cdot ch\,\,r)}{(z^{2} -1)^{2} sh^{2}\,\, x} .
\]

As, $ch\,\,(x-r)\le z\le ch\,\,(x+r),$ then the function $f(z)$ in point $z=ch\,\,x/ch\,\,r$ has maximum equal
\[
f_{\max } \left(\frac{ch\,\,x}{ch\,\,r} \right)=1-\left(\frac{ch^{2}\,\, x-ch^{2}\,\, r}{\sqrt{ch^{2}\,\, x-ch^{2}\,\, r} \cdot sh\,\,x} \right)^{2} =\]\[
=1-\frac{ch^{2}\,\, x-ch^{2}\,\, r}{sh^{2}\,\, x} =
\frac{ch^{2}\,\, r-1}{sh^{2}\,\, x} =\left(\frac{sh\,\,r}{sh\,\,x} \right)^{2} .
\]

Here from (27) we have
$$
A(x,r)\le \frac{1}{\lambda } \left(\frac{sh\,\,r}{sh\,\,x} \right)^{2\lambda } . \eqno(28)
$$%

Let $0<r\le c,$ then, making into account lemmas 1(а) and 2(а), for (25) with $0\le x\le r\le c$
$$
M_{G} f(ch\,\,x)=\mathop{\sup }\limits_{0<r\le 1} \frac{1}{\mu H(0,r)} \int\limits _{0}^{r}A_{ch\,\,t}^{\lambda } \left|f(ch\,\,x)\right|d\mu (t)
$$
$$
=\mathop{\sup }\limits_{0<r\le 1} \frac{\mu H(x,r)}{\mu H(0,r)} \cdot \frac{1}{\mu H(x,r)} \int\limits _{x-r}^{x+r}\left|f(ch\,\,t)\right|sh^{2\lambda }\,\,tdt
$$
$$
\le \mathop{c_{\lambda } \sup }\limits_{0<r\le 1} \frac{1}{\mu H(x,r)} \int\limits _{H(x,r)}\left|f(ch\,\,t)\right|d\mu (t)=c_{\lambda } M_{\mu } f(ch\,\,x).
\eqno(29)
$$%

And for $c<r<x<\infty $ from (28) and (20) we obtain 
\[
M_{G} f(ch\,\,x)\le \mathop{\sup }\limits_{0<r\le 1} \frac{A(x,r)\mu H(x,r)}{\mu H(0,r)\mu H(x,r)} \int\limits _{x-r}^{x+r}\left|f(ch\,\,t)\right|sh^{2\lambda }\,\,tdt
\]
\[
\le \mathop{c_{\lambda } \sup }\limits_{0<r\le 1} \frac{r\cdot ch^{2\lambda }\,\, x\cdot sh^{2\lambda }\,\, r}{\mu H(x,r)\left(sh\,\,\frac{r}{2} \right)^{2\lambda +1} sh^{2\lambda }\,\, x} \int\limits _{x-r}^{x+r}\left|f(ch\,\,t)\right|d\mu (t)
\]
\[
\le c_{\lambda } \left(\frac{ch\,\,x}{sh\,\,x} \right)^{2\lambda } \mathop{\sup }\limits_{0<r\le 1} ch^{2\lambda }\,\, \frac{r}{2} \cdot \frac{1}{\mu H(x,r)} \int\limits _{x-r}^{x+r}\left|f(ch\,\,t)\right|d\mu (t)
\]
\[
\le c_{\lambda } \left(\frac{e^{x} +e^{-x} }{e^{x} -e^{-x} } \right)^{2\lambda } ch^{2\lambda }\,\, \frac{1}{2} \sup \limits_{0<r\le 1} \frac{1}{\mu H(x,r)} \int\limits _{H(x,r)}\left|f(ch\,\,t)\right|d\mu (t)
\]
\[ \le c_{\lambda } \cdot 4^{\lambda } \cdot e\cdot M_{\mu } f(ch\,\,x),
\eqno\left(30\right)
\]
as $\dfrac{e^{2x} +1}{e^{2x} -1} \le 2\Leftrightarrow e^{2x} +1\le 2e^{2x} -2\Leftrightarrow e^{2x} \ge 3$ at $x\ge 1.$

From (29) and (30) follows, that
$$
M_{G} f(ch\,\,x)\le c_{\lambda } M_{\mu } f(ch\,\,x),0<r\le c,    \quad  0\le x<\infty . \eqno{(31)}
$$%

Now, we consider case, when $c<r<\infty.$

Point, that for $ch\,\,(x-r)\le z\le ch\,\,(x+r)$ and $x\ge 2r$ the function $f(z)=\dfrac{ch\,\,r-zch\,\,x}{\sqrt{z^{2} -1} sh\,\,x} $  has maximum equal $-\dfrac{\sqrt{ch^{2}\,\, x-ch^{2}\,\, r} }{sh\,\,x} $ .

Really,
$$
 f'(z)=-\frac{\sqrt{z^{2} -1} sh\,\,xch\,\,x+\frac{z}{\sqrt{z^{2} -1} } sh\,\,x(ch\,\,r-zch\,\,x)}{(z^{2} -1)sh^{2}\,\, x}
 $$
 $$
  =-\frac{(z^{2} -1)sh\,\,xch\,\,x+zsh\,\,xch\,\,r-z^{2} sh\,\,xch\,\,x}{(z^{2} -1)^{\frac{3}{2} } sh^{2}\,\, x}
=\frac{ch\,\,x-zch\,\,r}{(z^{2} -1)^{\frac{3}{2} } sh\,\,x} =0\Leftrightarrow z=\frac{ch\,\,x}{ch\,\,r} .
$$

In this point of the function $f(z)$ has maximum equal
$$
f_{\max }(z)=f \left(\frac{ch\,\,x}{ch\,\,r} \right)=\frac{ch^{2}\,\, r-ch^{2}\,\, x}{\sqrt{ch^{2}\,\, x-ch^{2}\,\, r} \cdot sh\,\,x} =-\frac{\sqrt{ch^{2}\,\, x-ch^{2}\,\, r} }{sh\,\,x}
$$
$$
=-\frac{ch\,\,x}{sh\,\,x} \sqrt{1-\left(\frac{ch\,\,r}{ch\,\,x} \right)^{2} } \sim -\frac{sh\,\,x}{ch\,\,x} , \eqno(32)
$$%
as
\[
\lim\limits_{x\to \infty } \frac{sh\,\,x}{ch\,\,x} =\lim\limits_{x\to \infty } \frac{e^{x} -e^{-x} }{e^{x} +e^{-x} } =1.
\]

From (27) and (32) we obtain
\[
 A(x,r)\le \int\limits _{-1}^{\frac{ch\,\,r-zch\,\,x}{\sqrt{z^{2} -1} \, \, \, sh\,\,x} }(1+u)^{\lambda -1} du\le \int\limits _{-1}^{-\frac{\sqrt{ch^{2}\,\, x-ch^{2} z} }{sh\,\,x} }(1+u)^{\lambda -1} du
 \]
 \[\sim \int\limits _{-1}^{-\frac{sh\,\,x}{ch\,\,x} }(1+u)^{\lambda -1} du=\frac{1}{\lambda } \left(1-\frac{sh\,\,x}{ch\,\,x} \right)^{\lambda } \le \frac{1}{\lambda } \left(1-\frac{sh^{2}\,\, x}{ch^{2}\,\, x} \right)^{\lambda }
 \]\[
 =\frac{1}{\lambda } \left(ch\,\,x\right)^{-2\lambda }  ,\ x\to\infty.
\eqno(33)
\]

Now, taking into account lemmas 1(b) and 2(b), also inequalities (24) and (33), we get
$$
\frac{\left|H(x,r)\right|_{\lambda } }{\left|H(0,r)\right|_{\lambda } } \le c_{\lambda } \left\{
\begin{array}{l} \dfrac{ch^{2\lambda }\,\, r}{ch^{4\lambda } \,\, \dfrac{r}{2} } , \\ \dfrac{ch^{2\lambda }\,\, xch^{2\lambda }\,\, r}{ch^{2\lambda }\,\, xch^{4\lambda } \,\, \dfrac{r}{2} }.  \end{array}\right. \le c_{\lambda },\quad c<r<\infty.
\eqno(34)
$$

Applying (34) we easy obtain
\[
M_{G} f(ch\,\,x)=\sup \limits_{r>c} \frac{1}{\left|H(0,r)\right|_{\lambda } } \int\limits _{0}^{r}A_{ch\,\,t}^{\lambda } \left|f(ch\,\,x)\right|d\mu (t)
\]
\[
=\mathop{\sup }\limits_{r>c} \frac{\left|H(x,r)\right|_{\lambda } }{\left|H(0,r)\right|_{\lambda } } \cdot \frac{1}{\left|H(x,r)\right|_{\lambda } } \int\limits _{\left|x-r\right|}^{x+r}\left|f(ch\,\,t)\right|sh^{2\lambda }\,\,tdt
 \]
 \[
 \le c_{\lambda } \frac{1}{\left|H(x,r)\right|_{\lambda } } \int\limits _{H(x,r)}\left|f(ch\,\,t)\right|d\mu  (t)=c_{\lambda } M_{\mu } f(ch\,\,x).
\eqno(35)
\]

Combine (31) and (35), we get
\[
M_{G} f(ch\,\,x)=\mathop{\sup }\limits_{0<r<\infty } \frac{1}{\left|H(0,r)\right|_{\lambda } } \int\limits _{0}^{r}A_{ch\,\,t}^{\lambda } \left|f(ch\,\,x)\right|d\mu (t)
\]\[
\le \sup\limits_{0<r\le c} \frac{1}{\left|H(0,r)\right|_{\lambda } } \int\limits _{0}^{r}A_{ch\,\,t}^{\lambda } \left|f(ch\,\,x)\right|d\mu (t)
 \]
 \[+\mathop{\sup }\limits_{r>c} \frac{1}{\left|H(0,r)\right|_{\lambda } } \int\limits _{0}^{r}A_{ch\,\,t}^{\lambda } \left|f(ch\,\,x)\right|d\mu (t)\le c_{\lambda } M_{\mu } f(ch\,\,x).
\]

Theorem 1 is proved.\\

\textbf{Theorem 2} {\it а) If $f\in L_{1,\lambda } \left[0,\infty \right),$ then for all  $\alpha >0$
$$
\left|\left\{x:M_{G} f(ch\,\,x)>\alpha \right\}\right|_{\lambda } \le \frac{c_{\lambda } }{\alpha } \int\limits _{0}^{\infty }\left|f(ch\,\,t)\right|sh^{2\lambda }\,\,tdt=\frac{c_{\lambda } }{\alpha } \left\| f\right\|  _{L_{1,\lambda }[0, \infty)} ,
$$
where $c_{\lambda } >0$ and depends only on $\lambda $.

b) If $f\in L_{p,\lambda } \left[0,\infty \right),$ $1<p<\infty ,$ then $M_{G} f(ch\,\,x)\in L_{p,\lambda } \left[0,\infty \right)$ and $\left\| M_{G} f\right\| _{L_{p,\lambda } \left[0,\infty \right)} \le c_{\lambda } \left\| f\right\| _{L_{p,\lambda } \left[0,\infty \right)} .$}

\textbf{Corollary 1} {\it If $f\in L_{p,\lambda } \left[0,\infty \right),$ $1\le p\le \infty ,$ then
$$
\lim \limits_{r\to 0} \frac{1}{\left|H(0,r)\right|_{\lambda } } \int\limits _{H(0,r)}A_{ch\,\,t}^{\lambda } f(ch\,\,x)sh^{2\lambda }\,\,tdt=f(ch\,\,x),
$$
for a. e. $x\in \left[0,\infty \right).$}\\

\textbf{Proof.} We need to introduce one maximal function defined on a space of homogeneous type. We mean a topological space $X$ equipped with a continuous preudometric $\rho$ and a positive measure  $\mu$, satisfying the doubling condition

\[\mu \left(E(x, 2r)\right)\leq C\mu \left(E(x, r)\right), \eqno(36) \]

with a constant  $C-$ independent of  $x$ and $r>0$.

Here $E(x,r)=\{y\in X:\rho(x,r)=|x-y|<r\}$.

Let $\left(X, \rho, \mu\right)$ is a space of homogeneous type. Let us define

\[M_{\mu}f(x)=\sup \limits_{r>0} \frac{1}{\mu E (x,r)} \int \limits_{E(x,r)} |f(t)|d\mu(t).\]

It is well known that the maximal function $ M_{\mu} $ is weak $ (1, 1) $ and is bounded on $ L_{p} \left(X, d\mu \right) $ for $ 1<p<\infty $ (see [19]). The measure of maximal function $ M_{\mu} f(ch\,\,x) $ introduced at the beginning of Section 1

\[\mu H(x,r)=|H(x,r)|_{\lambda}=\int \limits_{H(x,r)}sh^{2\lambda }\,\,t dt, \]

where

\[H(x,r)=\left\{\begin{array}{l} {(x-r, x+r) ,x-r>0;\, \, \ \ \ \, \, \, \, \, \, \, \, \, \, \, \, \,} \\ \\{(0, x+r),\quad x-r<0 ,\, \, \, \, \, \, \, \, \, \, \, \, \, \, \, \, \, \,}
\end{array}\right.\]
it is clear that this measure satisfies the condition (36), but then the confirmation of Theorem 2 follows from Theorem 1.

\textbf{The proof of the Corollary 1.} First show that for any function $ f\in L_{p,\lambda}[0, \infty), 1\leq p \leq \infty,$, representation $ ch\,\,t \mapsto A_{ch\,\,t}^{\lambda} f $ from $ \mathbb{R} $ into $ L_{p,\lambda} $ continuous, that is

\[\|A_{ch\,\,t}^{\lambda}f-f\|_{L_{p,\lambda}}\rightarrow 0 \,\,\,\,\mbox{при}\,\,\,\, t\rightarrow 0. \eqno {(37)}\]

Let $ f(x) $ is a continuous function definited for $ [a, b]\subset [0, \infty). $ Consider the function

\[ y(t, x, \varphi)=ch\,\,t ch\,\,x - sh\,\,tsh\,\,x cos\,\,\varphi \]

Hence we have

\[ |y(t, x, \varphi)  -  y(0, x, \varphi)|=|ch\,\,t ch\,\,x - sh\,\,t sh\,\,x cos\,\,\varphi -ch\,\,x|\]

\[=|(ch\,\,t-1) ch\,\,x - sh\,\,t sh\,\,x cos\,\,\varphi -ch\,\,x|\leq 2sh^{2}\,\,\frac{t}{2}ch\,\,x+2sh\,\,\frac{t}{2}ch\,\,\frac{t}{2}sh\,\,x\]

\[\leq 2sh\,\,\frac{t}{2}\left(sh\,\,\frac{t}{2}ch\,\,x+ch\,\,\frac{t}{2}sh\,\,x\right)=2sh\,\,\frac{t}{2}sh(\frac{t}{2}+x)\]

\[ \leq 2 sh\,\,\frac{t}{2}sh\left(\frac{t}{2}+b\right)\rightarrow 0 \,\,\,\,\mbox{при}\,\,\,\, t\rightarrow 0. \eqno {(38)}\]

On the strength of uniformly continuous of function $ f(x) $ on segment $ [a, b] $ for any $ \varepsilon >0 $ one may choose the number $ \delta >0 $, such that

\[|f[y(t, x, \varphi)]-f[y(0, x, \varphi)]| <\varepsilon,\,\,\,\mbox{have only}\]

\[|y(t, x, \varphi)-y(0, x, \varphi)|<\delta,\,\,\,\mbox{(that follows from (38)).}\]

Then we have

\[\left|A_{ch\,\,t}^{\lambda}f(ch\,\,x)-f(ch\,\,x)\right|\]

\[ \leq\frac{\Gamma\left(\lambda+\frac{1}{2}\right)}{\Gamma\left(\frac{1}{2}\right)\Gamma\left(\lambda\right)}\int_{0}^{\pi} |f[y(t,x, \varphi)]-f[y(0,x,\varphi)]|(sin\,\, \varphi)^{2\lambda-1}d\varphi<\varepsilon.\]

It follows, that

\[\|A_{ch\,\,t}^{\lambda}f-f\|_{\infty, \lambda}=\sup \limits_{x\in [a,b]}|A_{ch\,\,t}^{\lambda}f(ch\,\,x)-f(ch\,\,x)|<\varepsilon.\]

And for $ 1\leq p<\infty $

\[\|A_{ch\,\,t}^{\lambda}f-f\|_{L_{p, \lambda}[a,b]}=\left(\int_{a}^{b}|A_{ch\,\,t}^{\lambda}f(ch\,\,x)-f(ch\,\,x)|^{p}sh^{2\lambda}xdx\right)^{\frac{1}{p}}\]

\[<\varepsilon\left(\int_{a}^{b}sh^{2\lambda\,\,}xdx\right)^{\frac{1}{p}}<c_{p,\lambda}\varepsilon.\]

Thus for any continuous function definite by segment $ [a,b]\subset[0,\infty) $  and for any number $ \varepsilon >0 $ the inequality is valid:

\[\|A_{ch\,\,t}^{\lambda}f-f\|_{L_{p, \lambda}[a,b]}<\varepsilon\,\,\,\mbox{при}\,\,\, 1\leq p \leq \infty. \eqno{(39)}\]

It is known the set of all continuous functions with compact support in $ [0,\infty) $ is dense in $ L_{p,\lambda}[0,\infty) .$ Therefore for any number $ \varepsilon >0 $ there exists a continuous function with compact support in  $ [0,\infty) $, such that

\[\| f-f_{\varepsilon}\|_{ L_{p,\lambda}[0,\infty)}<\varepsilon. \eqno{(40)} \]

We denote $ g_{\varepsilon}=f-f_{\varepsilon}. $ Then $ g_{\varepsilon}\in L_{p,\lambda}[0,\infty)  $ and

\[\|g_{\varepsilon}\|_{L_{p,\lambda}[0,\infty)}<\varepsilon. \eqno{(41)}\]

Thus, if $f\in L_{p,\lambda}[0,\infty),$ then for any number $ \varepsilon >0 $ there exists a continuous function $ f_{\varepsilon} $ with the compact support and function $ g_{\varepsilon}\in L_{p,\lambda}[0,\infty)  $ with condition $ \|g_{\varepsilon}\|_{L_{p,\lambda}[0,\infty)}<\varepsilon,$ such that $ f=f_{\varepsilon}+g_{\varepsilon}. $

Hence we have $ A_{ch\,\,t}^{\lambda}f(ch\,\,x)=A_{ch\,\,t}^{\lambda}f_{\varepsilon}(ch\,\,x)+A_{ch\,\,t}^{\lambda}g_{\varepsilon}(ch\,\,x)-f(ch\,\,x)+f_{\varepsilon}(ch\,\,x)-f_{\varepsilon}(ch\,\,x), $ from which follows, that

\[\|A_{ch\,\,t}^{\lambda}f-f\|_{L_{p,\lambda}[0,\infty)}\leq\|A_{ch\,\,t}^{\lambda}f_{\varepsilon}-f_{\varepsilon}\|_{L_{p,\lambda}[0,\infty)}\]

\[+\|f-f_{\varepsilon}\|_{L_{p,\lambda}[0,\infty)}+\|A_{ch\,\,t}^{\lambda}g_{\varepsilon}\|_{L_{p,\lambda}[0,\infty)}.\]

Today, taking into account, that (see [17], lemma 2)

\[\|A_{ch\,\,t}^{\lambda}f_{\varepsilon}\|_{L_{p,\lambda}[0,\infty)}\leq\|f\|_{L_{p,\lambda}[0,\infty)}, t\in [0,\infty), 1\leq p \leq \infty\]
and also the inequality (39), (40) и (41), we get

\[\|A_{ch\,\,t}^{\lambda}f_{\varepsilon}-f\|_{L_{p,\lambda}[0,\infty)}\leq 3\varepsilon,\]
from which follows (37).

By virtue of the locality of the problem, one can account that $ f\in L_{1, \lambda}[0,\infty). $ In general case one can multiply $ f $ by characteristic function of ball $ B(0,r) $ and obtain required convergence almost everywhere interior to this ball and the tending $r$ to infinity could be got on ball interval $ [0,\infty). $

Suppose for any $ r>0 $ and for any $ x\in [0, \infty). $

\[ f_{r}(ch\,\,x)=\frac{1}{|B(0,r)|_{\lambda}} \int \limits_{B(0,r)}A_{ch\,\,t}^{\lambda}f(ch\,\,x)sh^{2\lambda}\,\,xdx.\]

Let $ r_{0}>0, B=B(0,r_{0}).$ According to Minkowski generalized inequality and discount (37), we obtain

\[\|f_{r}-f\|_{L_{1,\lambda}(B)}=\|\frac{1}{|B(0,r)|_{\lambda}} \int \limits_{B(0,r)}\left(A_{ch\,\,t}^{\lambda}f(ch\,\,x)-f(ch\,\,x)\right)sh^{2\lambda }\,\,tdt\|_{L_{1,\lambda}(B)} \]

\[\leq\frac{1}{|B(0,r)|_{\lambda}} \int \limits_{B(0,r)}\|A_{ch\,\,t}^{\lambda}f-f\|_{L_{1,\lambda}(B)}sh^{2\lambda }\,\,t dt \]

\[\leq \sup \limits_{|t|\leq r_{0}}\|A_{ch\,\,t}^{\lambda}f-f\|_{L_{1,\lambda}(B)}\rightarrow 0, \,\,\,\,\mbox{при}\,\,\,\,r_{o}\rightarrow + 0.\]

It means that there exists such sequence $ r_{k}, $ that $ r_{k}\rightarrow + 0, $ ($ k\rightarrow\infty  $) and $ \lim\limits_{k\rightarrow \infty} f_{r_k}(ch\,\,x)=f(ch\,\,x) $ almost everywhere in $ x\in [0, \infty).$

Now, let's prove that $ \lim\limits_{r\rightarrow + 0}f_{r}(ch\,\,x) $ exists almost everywhere. For this purpose for any $ x\in [0, \infty).$

\[\Omega _{f}(ch\,\,x)=\mid \overline{\lim\limits_{r\rightarrow +0}}f_{r}(ch\,\,x)-\lim\limits_{\overline{r\rightarrow + 0}}f_{r}(ch\,\,x)\mid\]

the oscillation of  $ {f_{r}} $ at the point $ x $ as $ r\rightarrow +0. $

If $ g  $ is a continuous function with compact support on $ [0,\infty) $, then $ g_r $ is convergent to $ g $  and consequently $ \Omega_{g}\equiv 0 $ is identically equal to zero in this case.

Further, if $ g\in L_{1,\lambda}[0, \infty) $, then, according to the statement of Theorem 2

\[|\left\{x\in [0,\infty):M_{G}g(ch\,\,x)>\varepsilon\right\}|_{\lambda}\leq\frac{c}{\varepsilon}\|g\|_{L_{1,\lambda}[0,\infty)},\,\, g\in L_{1,\lambda}[0,\infty).\]

On the other hand it is obvious that $ \Omega g(ch\,\,x)\leq 2M_{G}g(ch\,\,x).$ Thus

\[|\left\{x\in [0,\infty):\Omega_{g}(ch\,\,x)>\varepsilon\right\}|_{\lambda}\leq\frac{2c}{\varepsilon}\|g\|_{L_{1,\lambda}[0,\infty)},\,\, g\in L_{1,\lambda}[0,\infty).\]

How it was evidence above, any function $ f\in L_{p,\lambda}[0,\infty) $ can be written in form $ f=h+g, $ where $ h $ is continuous function and has compact support on $ [0, \infty) $, and $ g\in L_{p,\lambda}[0,\infty) $, moreover $ \|g\|_{L_{p,\lambda}[0,\infty)}<\varepsilon, $ for any $ \varepsilon>0. $ But $ \Omega\leq\Omega_{h}+\Omega_{g} $ и $ \Omega_{h}\equiv 0, $ however is continuous $ h. $ Therefore it follows that

\[|\left\{x\in [0,\infty):\Omega_{g}(ch\,\,x)>\varepsilon\right\}|_{\lambda}\leq\frac{c}{\varepsilon}\|g\|_{L_{1,\lambda}[0,\infty)}.\]

 Taking in inequality $ \|g\|_{L_{1,\lambda}[0,\infty)}<\varepsilon $ number $ \varepsilon $ chosen arbitrary small, we get $ \Omega f=0 $ almost everywhere on $ [0, \infty). $ Consequently, $ \lim\limits_{r\rightarrow 0} f_{r}(ch\,\,x) $ exists almost everywhere on  $ [0, \infty), $ what was confirmed.

 \textbf{Remark 1} Theorem 2 was proved earlier by W. C. Connett and A. L. Schwartz [7] for the Jacobi-type hypergroups.

 \textbf{Corollary 2} If $f\in L_{1,\lambda}[0,\infty)  $, then (see [18], Theorem 1)

 \[\lim\limits_{r\rightarrow 0}\frac{1}{\left(sh\,\, \frac{r}{2}\right)^{2\lambda+1}} \int_{0}^{r}\left|A_{ch\,\,t}^{\lambda}f(ch\,\,x)-f(ch\,\,x)\right|sh^{2\lambda }\,\,t dt =0,\]
 almost everywhere on $ x\in [0,\infty). $

 From here it follows that for any $ \varepsilon >0 $ find such $ \delta >0 $, that for all $ r<\delta $ the inequality is just:

 \[\frac{1}{\left(sh\,\, \frac{r}{2}\right)^{2\lambda+1}} \int_{0}^{r}\left|A_{ch\,\,t}^{\lambda}f(ch\,\,x)-f(ch\,\,x)\right|sh^{2\lambda }\,\,t dt <\varepsilon .\]

 But then discount of lemma 1 (a), we obtain

\[\left|\frac{1}{|H(0,r)|_{\lambda}} \int\limits_{H(0,r)}\left[A_{ch\,\,t}^{\lambda}f(ch\,\,x)-f(ch\,\,x)\right]sh^{2\lambda }\,\,t dt \right|\]

 \[\leq\frac{1}{\left(sh\,\, \frac{r}{2}\right)^{2\lambda+1}} \int_{0}^{r}\left|A_{ch\,\,t}^{\lambda}f(ch\,\,x)-f(ch\,\,x)\right|sh^{2\lambda }\,\,t dt <\varepsilon ,\]
for all $ r<\delta $, that means approval of Corollary 1.

\

\section{Some Morrey embeddings, associated with the Gegenbauer expansion}

We shall define function spaces, generated by the Gegenbauer expansion $ G $.

\textbf{Definition 1. [12]} {\it Let $1\le p<\infty ,\,\,0\le \gamma \le 2\lambda +1,$ $\left[r\right]_{1} =\min \left\{1,r\right\}.$ We denote by $L_{p,\lambda ,\gamma } (\left[0,\infty \right),G)$ Morrey-Gegenbauer spaces ($G$- Morrey spaces) and by $\widetilde{L}_{p,\lambda ,\gamma } (\left[0,\infty \right),G)$ modified $G$- Morrey spaces which are the sets of functions $f$ locally integrable on $\left[0,\infty \right)$ with finite norms}
\[
\begin{array}{l} {\left\| f\right\| _{L_{p,\lambda ,\gamma } (\left[0,\infty \right),G)} =\mathop{\sup }\limits_{x,r\in (0,\infty )} \left(\left(sh\,\,\frac{r}{2} \right)^{-\gamma } \int\limits _{H(0,r)}(A_{ch\,\,t}^{\lambda } \left|f(ch\,\,x)\right|)^{p} sh^{2\lambda }\,\,tdt \right)^{\frac{1}{p} } } \\ {\left\| f\right\| _{\widetilde{L}  _{p,\lambda ,\gamma (\left[0,\infty ),G\right)}} =\mathop{\sup }\limits_{x,r\in (0,\infty )} \left(\left[sh\,\,\frac{r}{2} \right]_{1}^{-\gamma } \int\limits _{H(0,r)}A_{ch\,\,t}^{\lambda } \left|f(ch\,\,x)\right|)^{p} sh^{2\lambda }\,\,tdt \right)^{\frac{1}{p} } .} \end{array}
\]

\textbf{Definition 2 [10]} {\it We denote by BMO $(\left[0,\infty \right),G)$ the BMO-Gegenbauer spaces ($G$-BMO space) as the set of functions locally integrable on $\left[0,\infty \right)$, with finite norm
$$
\left\| f\right\| _{*,G} =\mathop{\sup }\limits_{x,r\in \left[0,\infty \right)} \frac{1}{\left|H(0,r)\right|_{\lambda } } \int\limits _{H(0,r)}\left|A_{ch\,\,t}^{\lambda } f(ch\,\,x)-\mathop{f}\limits_{H(0,r)} (ch\,\,x)\right|sh^{2\lambda }\,\,tdt,
$$
where
\[
f_{H\left(0,r\right)} (ch\,\,x)=\frac{1}{\left|H(0,r)\right|_{\lambda } } \int\limits _{H(0,r)}A_{ch\,\,t}^{\lambda } \left|f(ch\,\,x)\right|sh^{2\lambda }\,\,tdt.
\]

Note that}
\[
\widetilde{L}_{p,\lambda ,0} \left[0,\infty \right)=L_{p,\lambda } \left[0,\infty \right),\,\,\, L_{p,\lambda ,2\lambda +1} \left[0,\infty \right)=L_{\infty ,\lambda .}
\]
$$
\widetilde{L}_{p,\lambda ,\gamma } \left[0,\infty \right)\subseteq L_{p,\lambda } \left[0,\infty \right)\ \mbox{and} \ \left\| f\right\| _{L_{p,\lambda } \left[0,\infty \right)} \le \left\| f\right\| _{\widetilde{L}_{p,\lambda ,\gamma } \left[0,\infty \right)} .
$$

\textbf{Lemma 3} {\it Let $1\le p<\infty ,$ $0\le \gamma \le 2\lambda +1$ and $\alpha p=2\lambda +1-\gamma $. Then}
$$
L_{p,\lambda,\gamma}[0,\infty)\subset L_{1,\lambda,2\lambda +1-\alpha }[0,\infty) ~~ \mbox{and} ~~ \left\| f\right\|_{L_{1,\lambda,2\lambda +1-\alpha }} \le c_{\lambda, p} \left\| f\right\| _{L_{p,\lambda,\gamma }}.
$$

\textbf{Proof}. Let $f\in L_{p,\lambda,\gamma}[0,\infty),$ $1\le p<\infty ,$ $0\le \gamma \le 2\lambda +1,$ ${1/p} +{1/q} =1$ and $\alpha p=2\lambda +1-\gamma .$

Applying Holder's inequality, we have
$$
\int\limits _{H(0,r)}A_{ch\,\,t}^{\lambda } \left|f(ch\,\,x)\right|sh^{2\lambda }\,\,tdt
$$$$
\le \left(\int\limits _{H(0,r)}(A_{ch\,\,t}^{\lambda } \left|f(ch\,\,x)\right| )^{p} sh^{2\lambda }\,\,tdt\right)^{\frac{1}{p} } \left(\int\limits _{H(0,r)}sh^{2\lambda }\,\,tdt \right)^{\frac{1}{q} }
.\eqno(42)
$$%

From Lemma 1 (b) it follows that for $r>c$
$$
\left|H(0,r)\right|_{\lambda } \le c_{\lambda } ch^{4\lambda } \,\, \frac{r}{2} <c_{\lambda } \left(ch\,\,\frac{r}{2} \right)^{2\lambda +1} \le c_{\lambda } \left(3sh\,\,\frac{r}{2} \right)^{2\lambda +1} =c_{\lambda } \left(sh\,\,\frac{r}{2} \right)^{2\lambda +1} . \eqno(43)
$$%

From Lemma 1 (а) and (43) it follows that for any $0<r<\infty $
$$
\left|H(0,r)\right|_{\lambda } \le c_{\lambda } \left(sh\,\,\frac{r}{2} \right)^{2\lambda +1} . \eqno(44)
$$%

Taking into account (44) and (42), we obtain
\[
\int\limits _{H(0,r)}A_{ch\,\,t}^{\lambda } \left|f(ch\,\,x)\right|sh^{2\lambda }\,\,tdt
\]\[
\le c_{\lambda ,p} \left(sh\,\,\frac{r}{2} \right)^{\frac{2\lambda +1}{q} } \left(\int\limits _{H(0,r)}(A_{ch\,\,t}^{\lambda } \left|f(ch\,\,x)\right|)^{p} sh^{2\lambda }\,\,tdt \right)^{\frac{1}{p} } .
\]

Further,
\[
\left(sh\,\,\frac{r}{2} \right)^{\alpha -2\lambda -1} \int\limits _{H(0,r)}A_{ch\,\,t}^{\lambda } \left|f(ch\,\,x)\right|sh^{2\lambda }\,\,tdt
\]\[
\le c_{\lambda ,p} \left(sh\,\,\frac{r}{2} \right)^{\alpha -2\lambda -1+\frac{2\lambda +1}{q} } \left(\int\limits _{H(0,r)}(A_{ch\,\,t}^{\lambda } \left|f(ch\,\,x)\right|)^{p} sh^{2\lambda }\,\,tdt \right)^{\frac{1}{p} }
\]
\[
\begin{array}{c} {=c_{\lambda ,p} \left(sh\,\,\frac{r}{2} \right)^{\alpha -2\lambda -1+(2\lambda +1)\left(1-\frac{1}{p} \right)} \left(\int\limits _{H(0,r)}(A_{ch\,\,t}^{\lambda } \left|f(ch\,\,x)\right|)^{p} sh^{2\lambda }\,\,tdt \right)^{\frac{1}{p} } } \\ {=c_{\lambda ,p} \left(sh\,\,\frac{r}{2} \right)^{\alpha -\frac{2\lambda +1}{p} } \left(\int\limits _{H(0,r)}(A_{ch\,\,t}^{\lambda } \left|f(ch\,\,x)\right|)^{p} sh^{2\lambda }\,\,tdt \right)^{\frac{1}{p} } } \\ {=c_{\lambda ,p} \left\{\left(sh\,\,\frac{r}{2} \right)^{-\gamma } \int\limits _{H(0,r)}(A_{ch\,\,t}^{\lambda } \left|f(ch\,\,x)\right|)^{p} sh^{2\lambda }\,\,tdt \right\}^{p} =c_{\lambda ,p} \left\| f\right\| _{L_{p,\lambda ,\gamma } } .} \end{array}
\]

Thus \[ f\in L_{1,\lambda, 2\lambda+1-\alpha}[0,\infty) ~~ \mbox{and} ~~ \|f\|_{L_{1,\lambda, 2\lambda+1-\alpha}} \leq c_{\lambda, p}\|f\|_{L_{p,\lambda, \gamma}}. \]

Lemma 3 is proved.

\section{RIESZ-GEGENBAUER POTENTIAL (($ R-G$) -POTENTIAL)}

In this Section the concept of potential of Riesz-Gegenbauer associated with the Gegenbauer differential operator $ G $ is introduced and its presentation of integrals is found. Moreover, for it the theorem of Sobolev type is proved. For the function $f,g\in L_{1,\lambda} \left[ 1,\infty \right)$ in [22] of Gegenbauer transformation is defined for function $P_{\gamma }^{\lambda} \left( x\right)$  and $Q_{\gamma }^{\lambda} \left( x\right)$  which are eigenfunctions of this operator $G$.
$$     	
F_{P} \left[ f\left( t\right) \right] \mapsto \hat{f} _{P} \left( \gamma
\right) =\int\limits_{1}^{\infty }f\left( t\right) P_{\gamma }^{\lambda}
\left( t\right) \left( t^{2} -1\right) ^{\lambda -\frac{1}{2} } dt,\eqno(45)
$$
$$ 	
F_{Q} \left[ f\left( t\right) \right] \mapsto \hat{f} _{Q} \left( \gamma
\right) =\int\limits_{1}^{\infty }f\left( t\right) Q_{\gamma }^{\lambda }
\left( t\right) \left( t^{2} -1\right) ^{\lambda -\frac{1}{2} } dt.\eqno(46)
$$

The inverses Gegenbauer transformations are defined by the formulas
$$
F_{P}^{-1} \left[ \hat{f} _{P} \left( \alpha \right) \right] \mapsto
f\left( x\right) =c^{\ast}_{\lambda } \int\limits_{1}^{\infty }\hat{f} _{P} \left(
\gamma \right) Q_{\gamma }^{\lambda} \left( x\right) \left( \gamma ^{2} -1\right)
^{\lambda -\frac{1}{2} } d\gamma ,\eqno(47)
$$
$$ 	
F^{-1}_{Q} \left[ \hat{f} _{Q} \left( \alpha \right) \right] \mapsto f\left(
x\right) =c_{\lambda } \int\limits_{1}^{\infty }\hat{f} _{Q} \left( \gamma
\right) P_{\gamma }^{\lambda} \left( x\right) \left( \gamma ^{2} -1\right)
^{\lambda -\frac{1}{2} } d\gamma ,\eqno(48)
$$
where

$c^{\ast}_{\lambda}= \dfrac{2^{\frac{3}{2} -\lambda } \sqrt{\pi }\Gamma(\lambda+1) \Gamma \left( \frac{%
1}{2} -\gamma \right) \Gamma \left( \frac{3+2\lambda }{4} \right) \left(
\Gamma \left( \lambda +\frac{1}{2} \right) \Gamma \left( \frac{5-2\lambda }{4%
} \right) \cos \pi \lambda \right) ^{-1} }{F\left( 1,\frac{1}{2} -\lambda ;%
\frac{5-2\lambda }{4} ;\frac{1}{2} \right) -F\left( 1,\frac{1}{2} -\lambda ;%
\frac{5-2\lambda }{4} ;\frac{1-2\lambda }{2} \right) } .$

Preliminary we prove one lemma

{\bf Lemma 4.} {\it Let  $f,g\in L_{1,\lambda } \left[ 1,\infty \right) \cap L_{2,\lambda } \left[
1,\infty \right)$. Then the equality is just}
$$
\int\limits_{1}^{\infty }f\left( x\right) A_{t}^{\lambda } g\left( x\right)
\left( x^{2} -1\right) ^{\lambda -\frac{1}{2} } dx =c^{\ast}_{\lambda }
\int\limits_{1}^{\infty }\hat{f} _{P} \left( \gamma \right) \left( \widehat{A_{t}^{\lambda } g} \right) _{P} \left( \gamma \right) \left(
\gamma ^{2} -1\right) ^{\lambda -\frac{1}{2} } d\gamma .
\eqno(49)
$$

{\bf Proof.} From (49) we have
$$
\int\limits_{1}^{\infty }f\left( x\right) A_{t}^{\lambda } g\left( x\right)
\left( x^{2} -1\right) ^{\lambda -\frac{1}{2} } dx=c^{\ast}_{\lambda }
\int\limits_{1}^{\infty }A_{t}^{\lambda } g\left( x\right) \left( x^{2}
-1\right) ^{\lambda -\frac{1}{2} } dx
$$
$$
\times \int\limits_{1}^{\infty }\hat{f} _{P} \left( \gamma \right) Q_{\gamma
}^{\lambda } \left( x\right) \left( \gamma ^{2} -1\right) ^{\lambda -\frac{1}{%
2} } d\gamma.
\eqno(56)
$$
Such (see the proof of Lemma 8 in [22])
$$
\int\limits_{1}^{\infty}\hat{f} _{P} \left( \gamma \right) Q_{\gamma }^{\lambda }
\left( x\right) \left( \gamma ^{2} -1\right) ^{\lambda -\frac{1}{2} }
d\gamma \lesssim \left\| f\right\| _{L_{2,\lambda }} ,
$$
then making the inequalite (see [22], Lemma 2)
$$
\left\| A_{t}^{\lambda } g\right\| _{L_{1,\lambda }} \leq \left\| g\right\|
_{L_{1,\lambda } },
$$
we obtain
$$
\left| \int\limits_{1}^{\infty }A_{t}^{\lambda } g\left( x\right) \left(
x^{2} -1\right) ^{\lambda -\frac{1}{2} } dx \int\limits_{1}^{\infty }\hat{f}
_{P} \left( \gamma \right) Q_{\gamma }^{\lambda } \left( x\right) \left(
\gamma ^{2} -1\right) ^{\lambda -\frac{1}{2} } d\gamma \right|
$$
$$
\leq \left\|
f\right\| _{L_{2, \lambda} } \int\limits_{1}^{\infty }\left| A_{t}^{\lambda }
g\left( x\right) \right| \left( x^{2} -1\right) ^{\lambda -\frac{1}{2} } dx
=\left\| f\right\| _{L_{2, \lambda} } \left\| A_{t}^{\lambda } g\right\|
_{L_{2, \lambda} } \leq \left\| f\right\| _{L_{2, \lambda}} \left\| g\right\|
_{L_{2, \lambda} } .
$$

In accord of theorem of Fubini
$$ 	
c^{\ast}_{^{\lambda } } \int\limits_{1}^{\infty }A_{t}^{\lambda } g\left( x\right)
\left( x^{2} -1\right) ^{\lambda -\frac{1}{2} } dx \int\limits_{1}^{\infty }%
\hat{f} _{P} \left( \gamma \right) Q_{\gamma }^{\lambda } \left( x\right)
\left( \gamma ^{2} -1\right) ^{\lambda -\frac{1}{2} } d\gamma
$$
$$
=c^{\ast}_{\lambda } \int\limits_{1}^{\infty }A_{t}^{\lambda } g\left( x\right)
Q_{\gamma }^{\lambda } \left( x\right) \left( x^{2} -1\right) ^{\lambda -%
\frac{1}{2} } dx \int\limits_{1}^{\infty }\hat{f} _{P} \left( \gamma \right)
\left( \gamma ^{2} -1\right) ^{\lambda -\frac{1}{2} } d\gamma
$$
$$
=c^{\ast}_{\lambda } \int\limits_{1}^{\infty }\left( \widehat{
A_{t}^{\lambda } g} \right) _{Q} \left( \gamma \right) \hat{f} _{P} \left(
\gamma \right) \left( \gamma ^{2} -1\right) ^{\lambda -\frac{1}{2} } d\gamma.
\eqno(51)
$$

Taking into account (51) in (50), we obtain (49).

Lemma 4 is proved.

{\bf Definition 3} {\it For $0<\alpha <2\lambda +1$  Riesz-Gegenbauer potential ($\left( R-G\right) $- potential)
$I_{G}^{\alpha } f\left(ch\,x\right)$  defined by the equality}
$$
I_{G}^{\alpha } f\left( ch\,\,x\right) =G^{-\frac{\alpha }{2} }
f\left( ch\,\,x\right). \eqno(52)
$$
Such  (see [8], p. 1933)
\begin{equation*}
G P_{\gamma }^{\lambda } \left( ch\,\,x\right) =\gamma \left( \gamma
+2\lambda \right) P_{\gamma }^{\lambda } \left( ch\,\,x\right) ,
\end{equation*}
then making selfadjoint of operator  $G$ (see [17], Lemma 4), we obtain for (45)
$$
\left( \widehat{G_{\lambda } f} \right) _{P} \left( \gamma \right)
=\int\limits_{1}^{\infty }P_{\gamma }^{\lambda } \left( ch\,\,x\right)
G f\left( ch\,\,x\right) sh^{2\lambda }\,\, xdx
$$$$
=\int\limits_{0}^{\infty
}f\left( ch\,\,x\right) \left( G_{\lambda } P_{\gamma }^{\lambda } \left(
ch\,\,x\right) \right) sh^{2\lambda }\,\, xdx
$$$$
=\gamma \left( \gamma +2\lambda \right) \int\limits_{\lambda }^{\infty
}f\left( ch\,\,x\right) P_{\gamma }^{\lambda } \left( ch\,\,x\right) sh^{2\lambda }
xdx =\gamma \left( \gamma +2\lambda \right) \hat{f}_{P} \left( \gamma
\right).
$$

Obviously, that by induction
$$
\left( \widehat{G_{\lambda }^{k} f} \right) _{P} \left( \gamma
\right) =\left( \gamma \left( \gamma +2\lambda \right) \right)^{k} \hat{f}
_{P} \left( \lambda \right) ,\qquad k=1,2, \ldots .
$$

This formula is naturally spread for the fractional indexes in the following form:
$$
\left( \widehat{G_{\lambda }^{-\frac{\alpha }{2} } f} \right) _{P}
\left( \gamma \right) :=\left( \gamma \left( \gamma +2\lambda \right)
\right) ^{-\frac{\alpha }{2} } \hat{f} _{P} \left( \lambda \right) .
\eqno(53)
$$
But then for (52) and (53) we have
$$
\left( \widehat{I_{G}^{\alpha } f} \right) _{P} \left( \gamma
\right) =\left( \gamma \left( \gamma +2\lambda \right) \right) ^{-\frac{%
\alpha }{2} } \hat{f}_{P} \left( \lambda \right).
\eqno(54)
$$

{\bf Lemma 5.} {\it Let $h_{r} \left( ch\,\,x\right)$  is the kernel associated with $G $  and \linebreak $0<\alpha <2\lambda +1$.
Then}
$$ 	
I_{G}^{\alpha } f( ch\,\,t) =\frac{1}{\Gamma (\frac{\alpha }{2}) }
\int\limits_0^\infty\left(\int\limits_{0}^{\infty } r^{\frac{\alpha }{2} -1} h_{r} \left(
ch\,\,x\right) dr\right) A_{ch\,\,t}^{\lambda } f\left( ch\,\,x\right) sh^{2\lambda }
xdx.\eqno(55)
$$

{\bf Proof.} Let
\begin{equation*}
\left( \hat{h} _{r} \right) _{Q} \left( \gamma \right) =e^{-\gamma \left(
\gamma +2\lambda \right) r} ,
\end{equation*}
then from (48) it follows, that
\begin{equation*}
h_{r} \left( ch\,\,x\right) =\int\limits_{1}^{\infty }e^{-\gamma \left( \gamma
+2\lambda \right) r} P_{\gamma }^{\lambda } \left( ch\,\,x\right) \left( \gamma
^{2} -1\right) ^{\lambda -\frac{1}{2} } d\gamma .
\end{equation*}

At by Lemma 4
\begin{equation*}
\int\limits_{0}^{\infty }h_{r} \left( ch\,\,x\right) A_{ch\,\,x}^{\lambda } f
\left( ch\,\,x\right) sh^{2\lambda }\,\, xdx
$$
$$
=c^{\ast}_{\lambda } \int\limits_{1}^{\infty
}e^{-\gamma \left( \gamma +2\lambda \right) r} \left(\widehat{ A_{ch\,\,t}^{\lambda }
f}\right) _{P} \left( \gamma \right) \left( \gamma ^{2} -1\right) ^{\lambda -%
\frac{1}{2} } d\gamma.
\end{equation*}

Thus we have
$$
\int\limits_{0}^{\infty }\int\limits_{0}^{\infty }r^{\frac{\alpha }{2} -1}
h_{r} \left( ch\,\,x\right) A_{ch\,\,t}^{\lambda } f\left( ch\,\,x\right)
sh^{2\lambda }\,\, xdxdr
$$
$$
=c^{\ast}_{\lambda } \int\limits_{1}^{\infty }\left(
\int\limits_{0}^{\infty }r^{\frac{\alpha }{2} -1} e^{-\gamma \left( \gamma
+2\lambda \right) r} dr\right) \left( \widehat{A_{ch\,\,t}^{\lambda
} f} \right) _{P} \left( \gamma \right) \left( \gamma ^{2} -1\right)
^{\lambda -\frac{1}{2} } d\gamma
$$
$$
=\left| \gamma \left( \gamma +2\lambda \right) r=t,dr=\frac{dt}{\gamma
\left( \gamma +2\lambda \right) } \right|
$$$$
 =c^{\ast}_{\lambda }
\int\limits_{1}^{\infty }\left( \int\limits_{0}^{\infty }e^{-t} t^{\frac{%
\alpha }{2} -1} dt\right) \left( \gamma \left( \gamma +2\lambda \right)
\right) ^{-\frac{\lambda }{2} } \left( \widehat{A_{ch\,\,t}^{\lambda } f} \right) _{P} \left( \gamma \right) \left( \gamma ^{2}
-1\right) ^{\lambda -\frac{1}{2} } d\gamma
$$
$$
=c^{\ast}_{\lambda } \Gamma \left( \frac{\alpha }{2} \right)
\int\limits_{1}^{\infty }\left( \gamma \left( \gamma +2\lambda \right)
\right) ^{-\frac{\alpha }{2} } \left( \widehat{A_{ch\,\,t}^{\lambda
} f} \right) _{P} \left( \gamma \right) \left( \gamma ^{2} -1\right)
^{\lambda -\frac{1}{2} } d\gamma .
$$

Taking into account that (see [23], lemma 2)
\begin{equation*}
\left( \widehat{A_{ch\,\,t}^{\lambda } f} \right) _{P} \left( \gamma
\right) =\hat{f} _{P} \left( \gamma \right) Q_{\gamma }^{\lambda } \left(
ch\, t\right)
\end{equation*}
for (54) and (47) we obtain
$$
\int\limits_{0}^{\infty }\int\limits_{0}^{\infty }r^{\frac{\alpha }{2} -1}
h_{r} \left( ch\,\,x\right) A_{ch\,\,t}^{\lambda } f\left(ch\,\,x\right)
sh^{2\lambda }\,\, xdxdr
$$$$
=c^{\ast}_{\lambda } \Gamma \left( \frac{\alpha }{2} \right)
\int\limits_{1}^{\infty }\left( \gamma \left( \gamma +2\lambda \right)
\right) ^{-\frac{\alpha }{2} } \hat{f} _{P} \left( \gamma \right) Q_{\gamma
}^{\lambda } \left( ch\,\,t\right) \left( \gamma ^{2} -1\right) ^{\lambda -\frac{%
\alpha }{2} } d\gamma
$$
$$
=\Gamma \left( \frac{\alpha }{2} \right) \int\limits_{0}^{\infty }\left(
\widehat{I_{G}^{\alpha } f} \right) _{P} \left( \gamma \right)
Q_{\gamma }^{\lambda } \left( ch\,\,t\right) \left( \gamma ^{2} -1\right)
^{\lambda -\frac{1}{2} } d\lambda =\Gamma \left( \frac{\alpha }{2} \right)
I_{G}^{\alpha } f\left( ch\,\,t\right) ,
$$
from it and for (47) it follows, that
\begin{equation*}
I_{G}^{\alpha } f\left( ch\,\,t\right) =\frac{1}{\Gamma \left( \frac{\alpha }{2}
\right) } \int\limits_{0}^{\infty }\left( \int\limits_{0}^{\infty }r^{\frac{%
\alpha }{2} -1} h_{r} \left( ch\,\,x\right) dr\right) A_{ch\,\,t}^{\lambda }
f\left( ch\,\,x\right) sh^{2\lambda }\,\, xdx.
\end{equation*}

Lemma 5 is proved.

\textbf{Corollary 2.} The following inequality is valid
$$
\left| I_{G}^{\alpha } f\left( ch\,\,t\right) \right| \lesssim
\int\limits_{0}^{\infty }\left| A_{ch\,\,t}^{\lambda } f\left( ch\,\,x\right)
\right| \left( sh\,\,x\right) ^{\alpha -2\lambda-1} sh^{2\lambda }\,\, xdx.
\eqno(56)
$$
Really, from formula (see [8], p. 1933)
\begin{equation*}
P_{\gamma }^{\lambda } \left( ch\,\,x\right) =\frac{\Gamma \left( \gamma
+2\lambda \right) \cos \pi \lambda}{\Gamma \left( \gamma \right) \Gamma
\left( \gamma +\lambda +1\right) } \left( 2ch\,\,x\right) ^{-\gamma -2\lambda }
$$
$$
\times F\left( \frac{\gamma }{2} +\lambda ,\frac{\gamma }{2} +\lambda +\frac{1}{2}
;\gamma +\lambda +1;\frac{1}{ch^{2}\,\, x} \right)
\end{equation*}
we have
\begin{equation*}
\left| P_{\gamma }^{\alpha } \left( ch\,\,x\right) \right| \lesssim \left(
ch\,\,x\right) ^{-\gamma -2\lambda }\,\,\,,
\end{equation*}
Since the function of Gauss $F\left( \alpha ,\beta;\gamma ;x\right)$  is convergence by appointed importance of parameters on the interval $\left[ 0,\infty\right) $,  (see [20], p. 1054).

Taking into account the last inequality, we estimate from above $h_{r} \left( ch\,\,x\right)$
$$
\left| h_{r} \left( ch\,\,x\right) \right| \lesssim \int\limits_{1}^{\infty
}e^{-\gamma \left( \gamma +2\lambda \right) r} \left( ch\,\,x\right) ^{-\gamma
-2\lambda } \left( \gamma ^{2} -1\right) ^{\lambda-\frac{1}{2} } d\gamma
$$$$
\lesssim
\int\limits_{0}^{\infty }e^{-\left( \gamma +1\right) \left( \gamma+1 +2\lambda
\right) r}\gamma^{\lambda-\frac 12}(ch\,\,x)^{-\gamma-2\lambda-1} d\gamma
$$$$
\lesssim e^{-r} \left( ch\,\,x\right) ^{-2\lambda -1} \int\limits_{0}^{\infty
}\gamma ^{\lambda -\frac{1}{2} } \left( ch\,\,x\right) ^{-\gamma } d\gamma
=\left| \frac{1}{ch\,\,x}\leq \frac{e}{e^{x+1} } \right|
$$
$$
\lesssim e^{-r} \left(
ch\,\,x\right) ^{-2\lambda -1} \int\limits_{0}^{\infty }e^{-\left( x+1\right)
\gamma } \gamma ^{\lambda-\frac{1}{2} } d\gamma
=\left| \left( x+1\right) \gamma =u\right|
$$$$
\lesssim e^{-r} \left( ch\,\,x\right)
^{-2\lambda -1} \int\limits_{0}^{\infty }e^{-u} u^{\lambda-\frac{1}{2} }
du=\Gamma \left( \lambda +\frac{1}{2} \right) e^{-r} \left( ch\,\,x\right)
^{-2\lambda -1} .
$$
Hence we have
$$
\int\limits_{0}^{\infty }r^{\frac{\alpha }{2} -1} h_{r} \left( ch\,\,x\right) dr
\lesssim \left( ch\,\,x\right) ^{-2\lambda -1} \int\limits_{0}^{\infty }r^{\frac{%
\alpha }{2} -1} e^{-r} dr
$$$$
=\Gamma \left( \frac{\alpha }{2} \right) \left(
ch\,\,x\right) ^{-2\lambda -1} \leq \Gamma \left( \frac{\alpha }{2} \right)
\left( ch\,\,x\right) ^{\alpha -2\lambda -1}
\leq \Gamma \left( \frac{\alpha }{2} \right) \left( sh\,\,x\right) ^{\alpha
-2\lambda -1} .
$$

Taking into account this inequality on (55), we obtain our approval.

\

\

\textbf{Sobolev type theorem for Riesz-Gegenbauer potential.}

We consider the Riesz-Gegenbauer fractional integral
\[
\Im _{G}^{\alpha } f(ch\,\,x)=\int\limits _{0}^{\infty }A_{ch\,\,t}^{\lambda } (sh\,\,x) ^{\alpha -2\lambda -1} f(ch\,\,t)sh^{2\lambda }\,\,tdt,\quad 0<\alpha <2\lambda +1.
\]

For $(R-G)$-potential the following analogue of Hardy-Littlewood-Sobolev theorem is valid.

\textbf{Theorem 3} {\it Let $0<\alpha <2\lambda +1,$ $1\le p<\dfrac{2\lambda +1}{\alpha } $ and $\dfrac{1}{p} -\dfrac{1}{q} =\dfrac{\alpha }{2\lambda +1} .$

а) If $f\in L_{p,\lambda } \left[0,\infty \right),$ then the integral $\Im _{G}^{\alpha } f$ is convergence absolutely for any $x\in \left[0,\infty \right).$

b) If $1<p<\dfrac{2\lambda +1}{\alpha } ,$ $f\in L_{p,\lambda } \left[0,\infty \right),$ then $\Im _{G}^{\alpha } \in L_{q,\lambda } \left[0,\infty \right)$ and $$
\|I^\alpha_G f\|_{L_{q,\lambda}[0,\infty)}\le\left\| \Im _{G}^{\alpha } f\right\| _{L_{q,\lambda } \left[0,\infty \right)} \le c_{\alpha ,\lambda ,p} \left\| f\right\| _{L_{p,\lambda } \left[0,\infty \right)},  \eqno(57)
$$
where $c_{\alpha ,\lambda ,p} -$positive constant, depending only on subscribed indexes.

с) If $f\in L_{1,\lambda } \left[0,\infty \right),$ $\dfrac{1}{q} =1-\dfrac{\alpha }{2\lambda +1} ,$ then}
$$
\left|\left\{t\in \left[0,\infty \right)\right. :\Im _{G}^{\alpha } f(ch\,\,t)>\left. \beta\right\}\right|_{\lambda } \le \left(\frac{c_{\alpha ,\lambda } }{\beta } \left\| f\right\| _{L_{1,\lambda } \left[0,\infty \right)} \right)^{q} ,\beta >0. \eqno(58)
$$

\textbf{Proof}. Let $f\in L_{p,\lambda } \left[0,\infty \right),$ $1\le p<\dfrac{2\lambda +1}{\alpha } ,$ $f_{1} (ch\,\,x)=f(ch\,\,x)\chi_{_{(0,1)}} (x),$
$$
f_{2} (ch\,\,x)=f(ch\,\,x)-f_{1} (ch\,\,x),\ \mbox{where} \ \chi_{_{(0,1)}} (x)=\left\{
\begin{array}{l} {1,} \\ {0,}
\end{array}\right. \ \ \begin{array}{l} {x\in (0,1),} \\ {x\in [1,\infty).}
\end{array}
$$

Then
$$
\Im _{G}^{\alpha } f(ch\,\,x)=\Im _{G}^{\alpha } f_{1} (ch\,\,x)+\Im _{G}^{\alpha } f_{2} (ch\,\,x)=\Im _{1} (ch\,\,x)+\Im _{2} (ch\,\,x).
$$

We estimate above $\Im _{1} (ch\,\,x).$
\[
\left|\Im _{1} (ch\,\,x)\right|\le \int\limits _{0}^{1}(sh\,\,x)^{\alpha -2\lambda -1} A_{ch\,\,t}^{\lambda } \left|f(ch\,\,x)\right|sh^{2\lambda }\,\,tdt
\]
\[
=\int\limits _{0}^{\infty }(sh\,\,x)^{\alpha -2\lambda -1}   \chi_{_{(0,1)}} (t)A_{ch\,\,t}^{\lambda } \left|f(ch\,\,x)\right|sh^{2\lambda }\,\,tdt.
\]

By Young inequality [23]
$$
\left\| \Im _{1} \left(ch\right)\left(\cdot \right)\right\| _{L_{p,\lambda } \left[0,\infty \right)} \le \left\| f(ch)(\cdot )\right\| _{L_{p,\lambda } \left[0,\infty \right)} \cdot \left\| \left|\cdot \right|^{\alpha -2\lambda -1} \chi_{_{(0,1)}} \right\| _{L_{1,\lambda } \left[0,\infty \right)} .
\eqno(59)
$$%

Here
$$
\left\| \left|\cdot \right|^{\alpha -2\lambda -1} \chi_{_{(0,1)}} \right\| _{L_{1,\lambda } } =\int\limits _{0}^{1}(sh\,\,t)^{\alpha -2\lambda -1} sh^{2\lambda }\,\,tdt
$$$$
\le \int\limits _{0}^{1}(sh\,\,t)^{\alpha -1} ch\,\,tdt=\int\limits _{0}^{1}(sh\,\,t)^{\alpha -1} d(sh\,\,t)=\frac{1}{\alpha } sh^{\alpha } 1.
\eqno(60)
$$%

From (58) and (59) it follows, that $\Im _{1} (ch\,\,x)$ for any $x\in \left[0,\infty \right)$ is convergence absolutely.

By using the Holder inequality
\[
\left|\Im _{2} \left(ch\,\,x\right)\right|\le \int\limits _{1}^{\infty }(sh\,\,t)^{\alpha -2\lambda -1} A_{ch\,\,t}^{\lambda } \left|f(ch\,\,x)\right|sh^{2\lambda }\,\,tdt
\]
\[
\le \left\| A_{ch\,\,t}^{\lambda } f\right\| _{L_{p,\lambda } } \cdot \left(\int\limits _{1}^{\infty }\left(sh\,\,t\right)^{\left(\alpha -2\lambda -1\right)q} sh^{2\lambda }\,\,tdt \right)^{\frac{1}{q} }
\]
\[
\le \left\| f\right\| _{L_{p,\lambda } } \left(\int\limits _{1}^{\infty }\left(sh\,\,t\right)^{\left(\alpha -2\lambda -1\right)q+2\lambda } ch\,\,tdt \right)^{\frac{1}{q} }
 \]\[
 =\left\| f\right\| _{L_{p,\lambda } } \left(\int\limits _{1}^{\infty }\left(sh\,\,t\right)^{\left(\alpha -2\lambda -1\right)q+2\lambda }  d(sh\,\,t)\right)^{\frac{1}{q} }
\]
\[
=\left(\frac{\left(sh1\right)^{\left(\alpha -2\lambda -1\right)q+2\lambda +1} }{\left(2\lambda +1-\alpha \right)q-2\lambda -1} \right)^{\frac{1}{q} } \cdot \left\| f\right\| _{L_{p,\lambda } } =c_{\alpha ,\lambda ,p} \left\| f\right\| _{L_{p,\lambda } } ,
\]
from here follows the absolutely convergence $\Im _{2} (ch\,\,x)$ for all $x\in \left[0,\infty \right).$

Thus, for all $f\in L_{p,\lambda } \left[0,\infty \right),$ $1\le p<\dfrac{2\lambda +1}{\alpha }\,\,\,  (R-G)$- potential $\Im _{G}^{\alpha } f(ch\,\,x)$ is convergence absolutely for all $x\in \left[0,\infty \right).$

b) We have
$$
\Im _{G}^{\alpha } f(ch\,\,x)=\left(\int\limits _{0}^{r}+\int\limits _{r}^{\infty }  \right)A_{ch\,\,t}^{\lambda } f(ch\,\,x)(sh\,\,t)^{\alpha -2\lambda -1} sh^{2\lambda }\,\,tdt
$$
$$
=A_{1} (x,r)+A_{2} (x,r). \eqno(61)
$$%

We consider $A_{1} (x,r).$
$$
 \left|A_{1} (x,r)\right|\le \int\limits _{0}^{r}\left|A_{ch\,\,t}^{\lambda } f(ch\,\,x)\right|(sh\,\,t)^{2\lambda } (sh\,\,t)^{\alpha -2\lambda -1} dt
$$
$$
  \le \sum _{k=0}^{\infty }\int\limits _{\frac{r}{2^{k+1} } }^{\frac{r}{2^{k} } }\frac{A_{ch\,\,t}^{\lambda } \left|f(ch\,\,x)\right|sh^{2\lambda }\,\,tdt}{(sh\,\,t)^{2\lambda +1-\alpha } }
 \le \sum _{k=0}^{\infty }\left(sh\frac{r}{2^{k+1} } \right)^{\alpha } \left(sh\frac{r}{2^{k+1} } \right)^{-2\lambda -1}
 $$
 $$
 \times \int\limits _{0}^{\frac{r}{2^{k} } }A_{ch\,\,t}^{\lambda } \left|f(ch\,\,x)\right|sh^{2\lambda }\,\,tdt \le \left(sh\,\,\frac{r}{2} \right)^{\alpha } M_{G} f(ch\,\,x).
\eqno(62)
$$%

We consider $A_{2} (x,r).$ By Holder inequality
\[
 \left|A_{2} (x,r)\right|\le \left(\int\limits _{r}^{\infty }\left|A_{ch\,\,t}^{\lambda } f(ch\,\,x)\right|^{p} sh^{2\lambda }\,\,tdt \right)^{\frac{1}{p} } \left(\int\limits _{r}^{\infty }(sh\,\,t)^{\left(\alpha -2\lambda -1\right)q} sh^{2\lambda }\,\,tdt \right)^{\frac{1}{q} }
\]
\[
  \le \left\| A_{ch\,\,t}^{\lambda } f\right\| _{L_{p,\lambda } } \left(\int\limits _{{r \mathord{\left/{\vphantom{r 2}}\right.\kern-\nulldelimiterspace} 2} }^{\infty }(sh\,\,t)^{\left(\alpha -2\lambda -1\right)q+2\lambda } ch\,\,tdt \right)^{\frac{1}{q} }
   \]
   \[
   \le \left\| f\right\| _{L_{p,\lambda } } \left(\frac{\left(sh\,\,\frac{r}{2} \right)^{\left(\alpha -2\lambda -1\right)q+2\lambda +1} }{\left(2\lambda +1-\alpha \right)q-2\lambda -1} \right)^{\frac{1}{q} }
=c_{\alpha ,\lambda ,p} \left\| f\right\| _{L_{p,\lambda } } \left(sh\,\,\frac{r}{2} \right)^{\alpha -2\lambda -1+\frac{2\lambda +1}{q} }
\]
\[
=c_{\alpha ,\lambda ,p} \left\| f\right\| _{L_{p,\lambda } } \left(sh\,\,\frac{r}{2} \right)^{\alpha -2\lambda -1+\left(2\lambda +1\right)\left(\frac{1}{p} -\frac{\alpha }{2\lambda +1} \right)}
 \]
 \[
 =c_{\alpha ,\lambda ,p} \left\| f\right\| _{L_{p,\lambda } } \left(sh\,\,\frac{r}{2} \right)^{\left(2\lambda +1\right)\left(\frac{1}{p} -1\right)} =c_{\alpha ,\lambda ,p} \left\| f\right\| _{L_{p,\lambda } } \left(sh\,\,\frac{r}{2} \right)^{-\frac{2\lambda +1}{q} } .
\eqno(63)
\]

Taking into account (62) and (63) in (61), we obtain
$$
\left|\Im _{G}^{\alpha } f(ch\,\,x)\right|\le c_{\alpha ,\lambda ,p} \left(\left(sh\,\, \frac{r}{2} \right)^{\alpha } M_{G} f(ch\,\,x)+\left(sh\,\,\frac{r}{2} \right)^{-\frac{2\lambda +1}{q} } \left\| f\right\| _{L_{p,\lambda } } \right).\eqno(64).
$$

Minimum of right-hand of the inequality (64) is reach for
$$
sh\,\,\frac{r}{2} =\left(\frac{\alpha q}{2\lambda +1} \cdot \frac{\left\| f\right\| _{L_{p,\lambda } } }{M_{G} f(ch\,\,x)} \right)^{\frac{p}{2\lambda +1} } .
$$

Then for (64) we have
\[
\left|\Im _{G}^{\alpha } f(ch\,\,x)\right|
\]\[
\le c_{\alpha ,\lambda ,p} \left\{\left(\frac{\left\| f\right\| _{L_{p,\lambda } } }{M_{G} f(ch\,\,x)} \right)^{\frac{\alpha p}{2\lambda +1} } M_{G} f(ch\,\,x)+\left(\frac{\left\| f\right\| _{L_{p,\lambda } } }{M_{G} f(ch\,\,x)} \right)^{-\frac{p}{q} } \left\| f\right\| _{L_{p,\lambda } } \right\}
\]
(as for the condition $\frac{1}{p} -\frac{1}{q} =\frac{\alpha }{2\lambda +1} \Rightarrow 1-\frac{p}{q} =\frac{\alpha p}{2\lambda +1} )=c_{\alpha ,\lambda ,p} \left(M_{G} f(ch\,\,x)\right)^{\frac{p}{q} } \left\| f\right\| _{L_{p,\lambda } }^{1-\frac{p}{q} } .$

From here we have
\[
\int\limits _{0}^{\infty }\left|\Im _{G}^{\alpha } f(ch\,\,t)\right|^{q} sh^{2\lambda }\,\,tdt\le c_{\alpha, \lambda, p} \left\|M_Gf(ch(\cdot)) \right\| _{L_{p,\lambda } }^{p} \cdot \left\| f\right\| _{L_{p,\lambda } }^{q-p}
\]
\[
\le c_{\alpha ,\lambda ,p} \left\| f\right\| _{L_{p,\lambda } }^{q-p} \cdot \left\| f\right\| _{L_{p,\lambda } }^{p} =c_{\alpha ,\lambda ,p} \left\| f\right\| _{L_{p,\lambda } }^{q} ,
\]
from here it follows
\[
\left\| \Im _{G}^{\alpha } f\right\| _{L_{p,\lambda } } \le c_{\alpha ,\lambda ,p} \left\| f\right\| _{L_{p,\lambda } } .
\]

с) Let $f\in L_{1,\lambda } \left[0,\infty \right).$ Denote
\[
\left|\left\{x\left. :\right|\Im _{G}^{\alpha } f(ch\,\,x)\right. \right|>2\left. \beta \right\}|_{\lambda } \le \left|\left\{x:\left|A_{1} (x,r)\right|>\beta \right\}|_{\lambda } +\left|\left\{x:\left|A_{2} (x,r)\right|>\beta \right.\right\} \right|_{\lambda } .
\]

From inequality (63) and Theorem 2 we have
\[
\beta \left|\left\{x\in \left[0,\infty \right)\left. :\left|A_{1} (x,r)\right.\right|>\beta \right\}\right|_{\lambda } =\beta \int\limits _{\left\{x\in \left[0,\infty \right)\left. :\left|A_{1} (x,r)\right.\right|>\beta \right\}}sh^{2\lambda }\,\, xdx
\]\[
\le \beta \int\limits _{\left\{x\in \left[0,\infty \right):c_{\alpha ,\lambda } \left(sh^{\alpha }\,\, \frac{r}{2} \right)M_{G} f(ch\,\,x)>\beta \right\}}sh^{2\lambda }\,\, xdx
 \]
 \[ =\beta \left|\left\{x\in \left[0,\infty \right):M_{G} f(ch\,\,x)>\frac{\beta }{c_{\alpha ,\lambda } sh^{\alpha }\,\, \frac{r}{2} } \right\}\right|_{\lambda }
 \]\[
 \le
 \beta \cdot \frac{c_{\alpha ,\lambda } }{\beta } sh^{\alpha }\,\, \frac{r}{2} \int\limits _{0}^{\infty }\left|f(ch\,\,x)\right|sh^{2\lambda }\,\, xdx
   =c_{\alpha ,\lambda } sh^{\alpha }\,\, \frac{r}{2} \left\| f\right\| _{L_{1,\lambda } } ,
\]
and also
\[
 \left|A_{2} (x,r)\right|\le \int\limits _{r}^{\infty }\left|A_{ch\,\,t}^{\lambda } f(ch\,\,x)\right|(sh\,\,t)^{\alpha -2\lambda -1} sh^{2\lambda }\,\,tdt
 \]
 \[
  \le \int\limits _{r}^{\infty }\frac{\left|A_{ch\,\,t}^{\lambda } f(ch\,\,x)\right|sh^{2\lambda }\,\,tdt}{\left(sh\,\,t\right)^{2\lambda +1-\alpha } } \le
   \int\limits _{r}^{\infty }\frac{\left|A_{ch\,\,t}^{\lambda } f(ch\,\,x)\right|sh^{2\lambda }\,\,tdt}{\left(sh\,\,\frac{t}{2} \right)^{2\lambda +1-\alpha } }
 \]
 \[
 \le \left(sh\,\,\frac{r}{2} \right)^{\alpha -2\lambda -1} \int\limits _{r}^{\infty }\left|A_{ch\,\,t}^{\lambda } f(ch\,\,x)\right|sh^{2\lambda }\,\,tdt \le \left(sh\,\,\frac{r}{2} \right)^{-\frac{2\lambda +1}{q} } \left\| f\right\| _{L_{1,\lambda } } .
\]

Suppose $\left(sh\,\,\frac{r}{2} \right)^{-\frac{2\lambda +1}{q} } \left\| f\right\| _{L_{1,\lambda } } =\beta ,$ we obtain $\left|A_{2} (x,r)\right|\le \beta $ and consequently $\left|\left\{x\in \left[0,\infty \right):\left|A_{2} (x,r)\right|>\beta \right\}\right|_{\lambda } =0.$

As last,
\[
\left|\left\{x\in \left[0,\infty \right):\left|\Im _{G}^{\alpha } f(ch\,\,x)\right|>\beta \right\}\right|_{\lambda } \le c_{\alpha ,\lambda } \cdot \frac{1}{\beta } sh^{\alpha }\,\,\, \frac{r}{2} \left\| f\right\| _{L_{1,\lambda } }
 \]\[
 \le c_{\alpha ,\lambda } \left(sh\,\,\frac{r}{2} \right)^{\alpha +\frac{2\lambda +1}{q} } =c_{\alpha ,\lambda } \left(sh\,\,\frac{r}{2} \right)^{2\lambda +1}
=c_{\alpha ,\lambda } \left(\frac{1}{\beta } \left\| f\right\| _{L_{1,\lambda } } \right)^{q} .
\]

Thus, $f\mapsto \Im _{G}^{\alpha } f$ is weak type $(1,q)$.

Theorem is proved,

\textbf{Theorem 4} {\it Let $0<\alpha <2\lambda +1,$ $p\alpha =2\lambda +1,$ $f\in L_{p,\lambda } \left[0,\infty \right),$ $\frac{1}{p} +\frac{1}{q} =1.$

Then $\widetilde{\Im}_{G}^{\alpha } f\in BMO\left[0,\infty \right)$ and the inequality
$$
\left\| \widetilde{\Im }_{G}^{\alpha } f\right\| _{BMO} \le c_{\alpha ,\lambda ,p} \left\| f\right\| _{L_{p,\lambda }} ,
$$
is fair, where $c_{\alpha ,\lambda ,p} >0-$ constant, depending only on written out indexes.}

\textbf{Proof}. We suppose
\[
f_{1} \left(ch\,\,x\right)=f\left(ch\,\,x\right)\chi_{_{(0,r/4)}} \left(ch\,\,x\right), f_{2} \left(ch\,\,x\right)=f\left(ch\,\,x\right)-f_{1} \left(ch\,\,x\right),
\]
where $\chi_{_{(0,r/4)}} \left(ch\,\,x\right)-$ is the characteristic function of the interval $\left[0,\infty \right),$ that is,
\[
\chi_{_{(0,r/4)}} \left(ch\,\,x\right)=\left\{\begin{array}{l} {1,} \\ {0,} \end{array}\right.          \begin{array}{l} {0\le x\le \frac{r}{4} ;} \\ {x>\frac{r}{4} .} \end{array}
\]

Then
\[
\widetilde{\Im}_{G}^{\alpha } f\left(ch\,\,x\right)=\widetilde{\Im }_{G}^{\alpha } f_{1} \left(ch\,\,x\right)+\widetilde{\Im }_{G}^{\alpha } f_{2} \left(ch\,\,x\right)=F_{1} \left(ch\,\,x\right)+F_{2} \left(ch\,\,x\right),
\]
where
\[
F_{1} \left(ch\,\,x\right)=\int\limits _{0}^{r/4 }\left(A_{ch\,\,t}^{\lambda } \left(sh\,\,x\right)^{\alpha -2\lambda -1} -\left(sh\,\,t\right)^{\alpha -2\lambda -1}
\chi_{_{\big(\frac{1}{4} ,\infty\big)}} \left(ch\,\,t\right)\right)f\left(ch\,\,t\right)sh^{2\lambda }\,\,tdt ,
\]
\[
F_{2} \left(ch\,\,x\right)=\int\limits _{r/4 }^{\infty }\left(A_{ch\,\,t}^{\lambda } \left(sh\,\,x\right)^{\alpha -2\lambda -1} -\left(sh\,\,t\right)^{\alpha -2\lambda -1} \chi_{_{\big(\frac{1}{4} ,\infty\big)}} \left(ch\,\,t\right)\right)f\left(ch\,\,t\right)sh^{2\lambda }\,\,tdt .
\]

As, function $f_{1} \left(ch\,\,x\right)$ has compact support, then the number
$$
a_{1} =-\int\limits _{(0,r/4)/\left(0,\min \left\{\frac{1}{4} ,\frac{r}{4} \right\}\right)}\left(sh\,\,t\right)^{\alpha -2\lambda -1} f\left(ch\,\,t\right)sh^{2\lambda }\,\,tdt
$$
is finite. And that is why we can write
$$
F_{1} \left(ch\,\,x\right)-a_{1} =\int\limits _{0}^{r/4 }A_{ch\,\,t}^{\lambda } \left(sh\,\,x\right)^{\alpha -2\lambda -1} f\left(ch\,\,x\right)sh^{2\lambda }\,\,tdt-
$$
$$
-\int\limits _{(0,r/4)/\left(0,\min \left\{\frac{1}{4} ,\frac{r}{4} \right\}\right)}\left(sh\,\,t\right)^{\alpha -2\lambda -1} f\left(ch\,\,t\right)sh^{2\lambda }\,\,tdt  +
$$
$$
+\int\limits _{(0,r/4)\left(0,\min \left\{\frac{1}{4} ,\frac{r}{4} \right\}\right. }\left(sh\,\,t\right)^{\alpha -2\lambda -1} f\left(ch\,\,t\right)sh^{2\lambda }\,\,tdt=
 $$$$
 =\int\limits _{0}^{r/4 }A_{ch\,\,t}^{\lambda } \left(sh\,\,x\right)^{\alpha -2\lambda -1} f\left(ch\,\,t\right)sh^{2\lambda }\,\,tdt=
 $$$$
  =\int\limits _{0}^{\infty }A_{ch\,\,t}^{\lambda } \left(sh\,\,x\right)^{\alpha -2\lambda -1} f_{1} \left(ch\,\,t\right)sh^{2\lambda }\,\,tdt.
\eqno(65)
$$%

Consider the integral
$$
A_{ch\,\,t}^{\lambda } f_{1} \left(ch\,\,x\right)=c_{\lambda } \int\limits _{0}^{\pi }\left|f\left(ch\,\,xch\,\,t-sh\,\,xsh\,\,t\cos\,\, \varphi  \right)\right|
$$
$$
\times \chi_{_{(0,r/4)}} \left(ch\,\,xch\,\,t-sh\,\,xsh\,\,t\cos\,\, \varphi  \right)\left(\sin\,\,\varphi \right)^{2\lambda -1} d\varphi  .
$$

So far as,
$$
ch\,\,\left(x-t\right)\le ch\,\,xch\,\,t-sh\,\,xsh\,\,t\cos\,\, \varphi  \le ch\,\,\left(x+t\right),
$$
then for $\left|x-t\right|>\frac{r}{4} $
\[
\chi_{_{(0,r/4)}} \left(ch\,\,xch\,\,t-sh\,\,xsh\,\,t\cos\,\, \varphi  \right)=0,
\]
and that is why
\begin{align*}
A_{ch\,\,t}^{\lambda } f_{1} \left(ch\,\,x\right)&=c_{\lambda } \int\limits _{\left\{\varphi \in \left[0,\pi \right],\left|x-t\right|\le r/4 \right\}}f\left(ch\,\,xch\,\,t-sh\,\,xsh\,\,t\cos\,\, \varphi  \right)\left(\sin\,\,\varphi \right)^{2\lambda -1} d\varphi
\\
& =A_{ch\,\,t}^{\lambda } f\left(ch\,\,x\right).
\end{align*}

But then for (65) we have
$$
\left|F_{1} \left(ch\,\,x\right)-a_{1} \right|\le \int\limits _{\left\{t\in \left[0,\infty \right):\left|x-t\right|\le r/4 \right\}}\left(sh\,\,t\right)^{\alpha -2\lambda -1} A_{ch\,\,t}^{\lambda } \left|f\left(ch\,\,x\right)\right|sh^{2\lambda }\,\,tdt . \eqno(66)
$$%

Study the estimation (66).

Let $\left(x-\frac{r}{4} ,x+\frac{r}{4} \right)\cap \left[0,\infty \right)=\left(0,x+\frac{r}{4} \right),$ then $0\le x\le r/4$ and we have for (66)
\[
 \left|F_{1} \left(ch\,\,x\right)-a_{1} \right|\le \int\limits _{0}^{x+r/4 }\left(sh\,\,t\right)^{\alpha -2\lambda -1} A_{ch\,\,t}^{\lambda } \left|f\left(ch\,\,x\right)\right|sh^{2\lambda }\,\,tdt
  \]\[
  \le \int\limits _{0}^{r}\left(sh\,\,t\right)^{\alpha -2\lambda -1} A_{ch\,\,t}^{\lambda } \left|f\left(ch\,\,x\right)\right|sh^{2\lambda }\,\,tdt
 \le c_{\alpha ,\lambda } \left(sh\,\,\frac{r}{2} \right)^{\alpha } M_{G} f\left(ch\,\,x\right).
\eqno(67)
\]

Let now $\left(x-\frac{r}{4} ,x+\frac{r}{4} \right)\cap \left[0,\infty \right)=\left(x-\frac{r}{4} ,x+\frac{r}{4} \right),$ then $x>\frac{r}{4} .$ Consider case, when $\frac{r}{4} \le x\le \frac{3r}{4} .$ Then
$$
\left|F_{1} \left(ch\,\,x\right)-a_{1} \right|\le \int\limits _{x-r/4 }^{x+r/4 }\left(sh\,\,t\right)^{\alpha -2\lambda -1} A_{ch\,\,t}^{\lambda } \left|f\left(ch\,\,x\right)\right|sh^{2\lambda }\,\,tdt
$$$$
\le \int\limits _{0}^{r}\left(sh\,\,t\right)^{\alpha -2\lambda -1} A_{ch\,\,t}^{\lambda } \left|f\left(ch\,\,x\right)\right|sh^{2\lambda }\,\,tdt
\le c_{\alpha ,\lambda } \left(sh\,\,\frac{r}{2} \right)^{\alpha } M_{G} f\left(ch\,\,x\right).
\eqno(68)
$$%

Finally, let, $\frac{3r}{4} \le x<\infty ,$ then by Holder inequality, we have
$$
\left|F_{1} \left(ch\,\,x\right)-a_{1} \right|=\int\limits _{x-\frac{r}{4} }^{x+\frac{r}{4} }\left(sh\,\,t\right)^{\alpha -2\lambda -1} A_{ch\,\,t}^{\lambda } \left|f\left(ch\,\,x\right)\right|sh^{2\lambda }\,\,tdt
$$$$
\le \left\| A_{ch\,\,t}^{\lambda } f\right\| _{L_{p,\lambda } } \left(\int\limits _{x-\frac{r}{4} }^{x+\frac{r}{4} }\left(sh\,\,t\right)^{\left(\alpha -2\lambda -1\right)q} sh^{2\lambda }\,\,tdt \right)^{\frac{1}{q} }
$$
$$
\le \left\| f\right\| _{L_{p,\lambda } } \left(sh\left(x-\frac{r}{4} \right)\right)^{\alpha -2\lambda -1} \left(\int\limits _{x-\frac{r}{4} }^{x+\frac{r}{4} }\left(2sh\,\,\frac{t}{2} ch\,\,\frac{t}{2} \right)^{2\lambda } dt \right)^{\frac{1}{q} }
 $$
 $$
 \le c_{\lambda ,p} \left\| f\right\| _{L_{p,\lambda } } \left(sh\left(x-\frac{r}{4} \right)\right)^{\alpha -2\lambda -1} \left(\int\limits _{x-\frac{r}{4} }^{x+\frac{r}{4} }\left(sh\,\,\frac{t}{2} \right)^{2\lambda } \left(ch\,\,\frac{t}{2} \right)^{2\lambda -1} d\left(sh\,\,\frac{t}{2} \right) \right)^{\frac{1}{q} }
$$
$$
\le c_{\lambda ,p} \left\| f\right\| _{L_{p,\lambda } } \left(sh\left(x-\frac{r}{4} \right)\right)^{\alpha -2\lambda -1} \left(\int\limits _{x-\frac{r}{4} }^{x+\frac{r}{4} }sh^{2\lambda }\,\, \frac{t}{2} d\left(sh\,\,\frac{t}{2} \right) \right)^{\frac{1}{q} }
 $$$$
 \le c_{\lambda ,p} \left\| f\right\| _{L_{p,\lambda } } \left(sh\left(x-\frac{r}{4} \right)\right)^{\alpha -2\lambda -1} \left(sh\left(\frac{x}{2} +\frac{r}{8} \right)\right)^{\frac{2\lambda +1}{q} }
$$
$$\le c_{\lambda ,p} \left\| f\right\| _{L_{p,\lambda } } \left(sh\left(\frac{x}{2} +\frac{r}{8} \right)\right)^{\alpha -2\lambda -1+\left(2\lambda +1\right)\left(1-\frac{\alpha }{2\lambda +1} \right)} =c_{\lambda ,p} \left\| f\right\| _{L_{p,\lambda } }.
\eqno(69)
$$

Combine (67), (68) and (69), we obtain that
\[
\left|F_{1} \left(ch\,\,x\right)-a_{1} \right|\le c_{\alpha ,\lambda } \left(sh\,\,\frac{r}{2} \right)^{\alpha } M_{G} f\left(ch\,\,x\right)+c_{\lambda ,p} \left\| f\right\| _{L_{p,\lambda } } .
\]

From here it follows, that
$$
\mathop{\sup }\limits_{r>0} \frac{1}{\left|H\left(0,r\right)\right|_{\lambda } } \int\limits _{0}^{r}\left|A_{ch\,\,t}^{\lambda } \left(F_{1} \left(ch\,\,x\right)-a_{1} \right)\right|sh^{2\lambda }\,\,tdt
 $$
 $$
 \le \mathop{\sup }\limits_{r>0} \frac{1}{\left|H\left(0,r\right)\right|_{\lambda } } \int\limits _{0}^{r}A_{ch\,\,t}^{\lambda } \left|F_{1} \left(ch\,\,x\right)-a_{1} \right|sh^{2\lambda }\,\,tdt
$$
$$
\le c_{\alpha ,\lambda } \sup\limits_{r>0} \left(sh\,\,\frac{r}{2} \right)^{\alpha -2\lambda -1} \int\limits _{0}^{r}A_{ch\,\,t}^{\lambda } \left|M_{G} f\left(ch\,\,x\right)\right|sh^{2\lambda }\,\,tdt
 $$
 $$
 +c_{\lambda ,p} \frac{\left\| f\right\| _{L_{p,\lambda } } }{\left|H\left(0,r\right)\right|_{\lambda } } \int\limits _{0}^{r}sh^{2\lambda }\,\,tdt
\le c_{\alpha ,\lambda } \mathop{\sup }\limits_{r>0} \left(sh\,\,\frac{r}{2} \right)^{\alpha -2\lambda -1} \left(\int\limits _{0}^{r}sh^{2\lambda }\,\,tdt \right)^{\frac{1}{q} }
$$
$$
\times\left\| A_{ch\,\,t}^{\lambda } \left|M_{G} f\left(\cdot \right)\right|\right\| _{L_{p,\lambda } } +c_{\lambda ,p} \left\| f\right\| _{L_{p,\lambda } }
$$
$$
\le c_{\alpha ,\lambda } \mathop{\sup }\limits_{r>0} \left(sh\,\,\frac{r}{2} \right)^{\alpha -2\lambda -1} \left(sh\,\,\frac{r}{2} \right)^{\frac{2\lambda +1}{q} } \left\| M_{G} f\right\| _{L_{p,\lambda } }+c_{\lambda ,p} \left\| f\right\| _{L_{p,\lambda } } \le c_{\alpha ,\lambda ,p} \left\| f\right\| _{L_{p,\lambda } } . \eqno(70)
$$

We suppose
$$
a_{2} =\int\limits _{\left(0,\max \left\{\frac{1}{4} ,\frac{r}{4} \right\}\right)/\left(0,\frac{r}{4} \right)}\left(sh\,\,t\right)^{\alpha -2\lambda -1} f\left(ch\,\,t\right)sh^{2\lambda }\,\,tdt .
$$

We estimate above the difference
$$
\left|F_{2} \left(ch\,\,x\right)-a_{2} \right|
$$$$
=\left|\int\limits _{\frac{r}{4} }^{\infty }\left(A_{ch\,\,t}^{\lambda } \left(sh\,\,x\right)^{\alpha -2\lambda -1} -\left(sh\,\,t\right)^{\alpha -2\lambda -1} \aleph _{\left(\frac{1}{4} ,1\right)} \left(ch\,\,t\right)\right)f\left(ch\,\,t\right)sh^{2\lambda }\,\,tdt \right.
 $$$$
 -\int\limits _{(0,\max \left\{\frac{1}{4} ,\frac{r}{4} \right\})/\left(0,\frac{r}{4} \right)}\left(sh\,\,t\right)^{\alpha -2\lambda -1} f\left(ch\,\,t\right)sh^{2\lambda }\,\,td \left. t\right|
 $$$$
 =\left|\int\limits _{\frac{r}{4} }^{\infty }\left(A_{ch\,\,t}^{\lambda } \left(sh\,\,x\right)^{\alpha -2\lambda -1} -\left(sh\,\,t\right)^{\alpha -2\lambda -1} \right)f\left(ch\,\,t\right)sh^{2\lambda }\,\,tdt \right|
 $$
 $$
 \le \int\limits _{\frac{r}{4} }^{\infty }\left|f\left(ch\,\,t\right)\right|B\left(x,t\right)sh^{2\lambda }\,\,tdt =\Im \left(x,r\right). \eqno(71)
$$

We consider expansion
$$
 B\left(x,t\right)=\left|A_{ch\,\,t}^{\lambda } \left(sh\,\,x\right)^{\alpha -2\lambda -1} -(sh\,\,t)^{\alpha-2\lambda-1}\right|
$$$$
 =c_{\lambda } \left|\int\limits _{0}^{\pi }\left(\left(ch\,\,xch\,\,t-sh\,\,xsh\,\,t\cos\,\, \varphi  \right)^{2} -1\right)^{\frac{\alpha -2\lambda -1}{2} } -\left(sh\,\,t\right)^{\alpha -2\lambda -1} ) \left(\sin\,\,\varphi \right)^{2\lambda -1} d\varphi \right|
$$$$
\le c_{\lambda } \int\limits _{0}^{\pi }\left|\left(\max \left(sh\left(x+t\right),\left|sh\left(x-t\right)\right|\right)\right)^{\alpha -2\lambda -1} -\left(sh\,\,t\right)^{\alpha -2\lambda -1} \right| \left(\sin\,\,\varphi \right)^{2\lambda -1} d\varphi .
$$

We estimate above the value $B(x,t)$. Easy to notice, that
$$
B(x,t)\lesssim \left|max\left(\left\{sh\,\,(x+t),|sh\,\,(x-t)|\right\}\right)^{\alpha-2\lambda-1}-(sh\,\,t)^{\alpha-2\lambda-1}\right|\equiv V(x,t).
\eqno(72)
$$

I	Let  $0<t<x-t<\infty $, then  $0<t<\dfrac{x}{2} <x+t$.

From here it follows, that
$$
\left( sh\,\,t\right) ^{\alpha -2\lambda -1} >\left( sh\left( x+t\right) \right)
^{\alpha -2\lambda -1 }.\eqno(73)
$$

II	Let  $0<x-t<t<\infty $, then $\dfrac{x}{2} <t<x<x+t,$ and in this case the inequality (73) is just.

III	Let   $0<t-x<\infty $, then  $x<t<x+t<\infty $.

Again the inequality (73) takes place.

IV	Let  $0<x+t<\infty $, as  $t<x+t$, then and here (73) is fair.

Combine all these cases, we obtain, that
\begin{equation*}
V\left( x,t\right) =\left( sh\,\,t\right) ^{\alpha -2\lambda -1} -\left(
sh\left( x+t\right) \right) ^{\alpha -2\lambda -1}.
\end{equation*}

Applying the Lagrange formula to segment $\left[ t,x+t\right] $, we obtain
\begin{equation*}
V\left( x,t\right) \equiv V_{\xi } \left( x,t\right) =\frac{\left( 2\lambda
+1-\alpha \right) xch\,\,\xi }{\left( sh\,\,\xi \right) ^{2\lambda +2-\alpha } }
,\quad t<\xi <t+x.
\end{equation*}

From here we have
\begin{equation*}
\nu _{\xi } \left( x,t\right) \lesssim \left\{
\begin{array}{l}
x\left( sh\,\,t\right) ^{\alpha -2\lambda -2}, \ \mbox{if}\ \xi <1,\qquad\qquad\qquad(74)
\\
x\left( sh\,\,t\right) ^{\alpha -2\lambda -1}, \ \mbox{if}\  \xi \geq 1\,\qquad\qquad\qquad\ (75)%
\end{array}
\right. \;\,
\end{equation*}

At first we consider the case, when  $\xi <1$.

Applying the Holder inequality and also(74) and (72), from (71) for  $x\leq r$ we obtain.
$$
\tau \left( x,r\right) =\int\limits_{r/4}^{\infty }\left| f\left( ch\,\,t\right)
\right| B\left( x,t\right) sh^{2\lambda }\,\,tdt\lesssim \left\| f\right\|
_{L_{p,\lambda } } x\left( \int\limits_{r/4}^{\infty }\frac{sh^{2\lambda }\,\,\
tdt}{\left( sh\,\,t\right) ^{\left( 2\lambda +2-\alpha \right) q} } \right) ^{%
\frac{1}{q} }
$$
$$
\lesssim \left\| f\right\| _{L_{p,\lambda } } r\left( \int\limits_{r/4}^{\infty
}\left( sh\,\,t\right) ^{\left( \alpha -2\lambda -2\right) q+2\lambda } dsh\,\,t
\right) ^{\frac{1}{q} } =\left\| f\right\| _{L_{p.\lambda } } \frac{r}{sh\,\,\
\frac{r}{4} } \lesssim \left\| f\right\| _{L_{p,\lambda } } .\eqno(76)
$$

As according to condition of theorem
$\alpha -2\lambda -2+\left( 2\lambda +1\right) /q=\alpha -2\lambda -2\linebreak +\left(
2\lambda +1\right) \left( 1-\dfrac{\alpha }{2\lambda +1} \right) =\alpha
-2\lambda -2+2\lambda +1-\alpha =-1.
$

 Now we consider the case, when  $\xi \geq 1$.

Let at first  $0<x\leq 8$. Taking into account (75) and acting how above for $x\leq r$ we obtain
$$
\tau \left( x,r\right) =\int\limits_{r/4}^{\infty }\left| f\left( ch\,\,t\right)
\right| B\left( x,t\right) sh^{2\lambda }\,\,tdt\lesssim \left\| f\right\|
_{L_{p,\lambda } } x\left( \int\limits_{r/4}^{\infty }\frac{sh^{2\lambda }\,\
tdt}{\left( sh\,\,t\right) ^{\left( 2\lambda +1-\alpha \right) q} } \right) ^{%
\frac{1}{q} }
$$
$$
\lesssim \left\| f\right\| _{L_{p,\lambda } } x\left( \int\limits_{x/4}^{\infty }%
\frac{\left( sh^{2\lambda }\,\, \frac{t}{2} \right) \left( ch\,\,\frac{t}{2} \right)
^{2\lambda -1} d(sh\,\,\frac{t}{2} )}{\left( 2sh\,\,\frac{t}{2} ch\,\,\frac{t}{2} \right)
^{\left( 2\lambda +1-\alpha \right) q} } \right) ^{\frac{1}{q} }
\lesssim
\left\| f\right\| _{L_{p,\lambda } } x\left( \int\limits_{x/4}^{\infty }%
\frac{\left( sh\,\,\frac t2\right)^{4\lambda -1} d (sh\,\,\frac{t}{2} )}{\left( sh\,\,\frac{t}{2}
\right) ^{\left( 2\lambda +1-\alpha \right) q} } \right) ^{\frac{1}{q} }
$$
$$
=\left\| f\right\| x\left( \int\limits_{x/4}^{\infty }\left( sh\,\,\frac{t}{2}
\right) ^{\left( \alpha -2\lambda -1\right) q+4\lambda -1} d(sh\,\,\frac{t}{2})
\right) ^{\frac{1}{q} }\lesssim \left\| f\right\| _{L_{p,\lambda } } x\left( sh\,\,\
\frac{x}{8} \right) ^{\alpha -2\lambda -1+4\lambda /q}
$$
$$
=\left\| f\right\| _{L_{p,\lambda } } x\left( sh\,\,\frac{x}{8} \right) ^{\alpha
-2\lambda -1+4\lambda  \left( 1-\frac{\alpha }{2\lambda +1} \right)}
=\left\| f\right\| _{L_{p,\lambda } } x\left( sh\,\,\frac{t}{8} \right) ^{\frac{%
1-2\lambda }{1+2\lambda } \alpha +2\lambda -1} $$
$$
\lesssim \left\| f\right\| _{L_{p,\lambda } } \left( sh\,\,\frac{x}{8} \right) {^{%
\frac{1-2\lambda }{1+2\lambda } \alpha +2\lambda } } \lesssim \left\| f\right\|
_{L_{p,\lambda } }. \eqno(77)
$$

Let now $8<x<\infty $. Then, how above, we obtain
$$
\tau \left( x,r\right) \lesssim \left\| f\right\| _{L_{p,\lambda } } x\left(
\int\limits_{x/4}^{\infty }\frac{sh^{2\lambda }\,\,tdt}{\left( sh\,\,t\right)
^{\left( 2\lambda +1-\alpha \right) q} } \right) ^{\frac{1}{q} } \lesssim
\left\| f\right\| _{L_{p,\lambda } } x\left( sh\,\,\frac{x}{8} \right) ^{\frac{%
1-2\lambda }{1+2\lambda } \alpha +2\lambda -1}
$$
$$
=\left\| f\right\| _{L_{p,\lambda } } \frac{x}{\left( sh\,\,\frac{x}{8} \right)
^{1-2\lambda -{\frac{1-2\lambda }{1+2\lambda } \alpha } } } =\left\|
f\right\| _{L_{p,\lambda } } \frac{x}{\left(2 sh\,\,\frac{x}{16} ch\,\,\frac{x}{16}
\right) ^{1-2\lambda -{\frac{1-2\lambda }{1+2\lambda } \alpha } } }
$$
$$
=\left\| f\right\| _{L_{p,\lambda } } \frac{x}{\left( 2sh\,\,\frac{x}{2^{n} }
\right) ^{2\left( 1-2\lambda -\frac{1-2\lambda }{1+2\lambda } x\right) } }
\lesssim ...\lesssim \left\| f\right\| _{L_{p,\lambda } } \frac{x}{\left( 2^{n} sh\,\,\
\frac{x}{2^{n+3} } \right) ^{2^{n} \left( 1-2\lambda -\frac{1-2\lambda }{%
1+2\lambda } x\right) } }
$$
$$
\lesssim \left\| f\right\| _{L_{p,\lambda } } \frac{x}{\left( \frac{x}{8}
\right) ^{2^{n} \left( 1-2\lambda -\frac{1-2\lambda }{1+2\lambda } \alpha
\right) } } \!\lesssim\! \left\| f\right\| _{L_{p,\lambda } } \frac{8}{\left( \frac{x%
}{8} \right) ^{2^{n} \left( 1-2\lambda -\frac{1-2\lambda }{1+2\lambda }
\alpha \right) -1} } \!\lesssim\! \left\| f\right\| _{L_{p,\lambda } } ,
\eqno(78)
$$

so far as for anough greater $n=n_{0} $
\begin{equation*}
2^{n_{0} } \left( 1-2\lambda -\frac{1-2\lambda }{1+2\lambda } \alpha \right)
-1\geq 0\Leftrightarrow \frac{1-2\lambda }{1+2\lambda } \alpha \leq
1-2\lambda -\frac{1}{2^{n_{0} } }
\end{equation*}
\begin{equation*}
\frac{1-2\lambda }{1+2\lambda } \alpha <1-2\lambda\Leftrightarrow \alpha <2\lambda
+1.
\end{equation*}

Combine the estimates (76), (77) and (78), for \linebreak $0<x\leq r$ on (71) we obtain
$$
\left| F_{2} \left( ch\,\,x\right) -a_{2} \right| \lesssim \left\| f\right\|
_{L_{p,\lambda } } .
$$

Hence we have
\begin{equation*}
\left| A_{ch\,\,t}^{\lambda } F_{2} \left( ch\,\,x\right) -a_{2} \right| \leq A_{ch\,\,t}^{\lambda } \left| F_{2} \left( ch\,\,x\right) -a_{2} \right| \lesssim \left\|
f\right\| _{L_{p,\lambda } } .\eqno(79)
\end{equation*}

From (79) it follows, that
$$
\sup \limits_{r>0}\frac{1}{\left| H\left( 0,r\right) \right|_\lambda }
\int\limits_{0}^{r}\left| A_{ch\,\,t}^{\lambda } F_{2} \left( ch\,\,x-a_{2} \right)
\right| sh^{2\lambda }\,\,tdt
$$
$$
\leq \sup \limits_{r>0}\frac{1}{\left| H\left(
0,r\right) \right| _{\lambda } }
\int\limits_{0}^{r}A_{ch\,\,t}^{\lambda } \left| F_{2} \left( ch\,\,x\right) -a_{2}
\right| sh^{2\lambda }\,\,tdt
$$
$$
\lesssim \left\|f\right\| _{L_{p,\lambda}} \sup\limits_{r>0} \frac{1}{\left| H\left( 0,r\right) \right|_\lambda } \int\limits_{0}^{r}sh^{2%
\lambda }\, tdt \lesssim \left\| f\right\|_{L_{p,\lambda}}.\eqno (80)
$$

Denote by
$$
a_{f} =a_{1} +a_{2} =\int\limits _{\left(0,\max \left\{\frac{1}{4} ,\frac{r}{4} \right\}\right)}\left(sh\,\,t\right)^{\alpha -2\lambda -1}  f\left(ch\,\,t\right)sh^{2\lambda }\,\,tdt.
$$

At last, from (70) and (80) finally we obtain
\[
\mathop{\sup }\limits_{r>0} \frac{1}{\left|H\left(0,r\right)\right|_{\lambda } } \int\limits _{0}^{r}\left|A_{ch\,\,t}^{\lambda } \tilde{\Im }_{G}^{\alpha } f\left(ch\,\,x\right)-a_{f} \right|sh^{2\lambda }\,\,tdt
\]
\[
 =\mathop{\sup }\limits_{r>0} \frac{1}{\left|H\left(0,r\right)\right|_{\lambda } } \int\limits _{0}^{r}\left|A_{ch\,\,t}^{\lambda } F_{1} \left(ch\,\,x\right)-a_{1} +A_{ch\,\,t}^{\lambda } F_{2} \left(ch\,\,x\right)-a_{2} \right| sh^{2\lambda }\,\,tdt
 \]
 \[
  \le \mathop{\sup }\limits_{r>0} \frac{1}{\left|H\left(0,r\right)\right|_{\lambda } } \int\limits _{0}^{r}\left|A_{ch\,\,t}^{\lambda } F_{1} \left(ch\,\,x\right)-a_{1} \right| sh^{2\lambda }\,\,tdt
  \]\[
  +\mathop{\sup }\limits_{r>0} \frac{1}{\left|H\left(0,r\right)\right|_{\lambda } } \int\limits _{0}^{r}\left|A_{ch\,\,t}^{\lambda } F_{2} \left(ch\,\,x\right)-a_{2} \right|sh^{2\lambda }\,\,tdt\lesssim \left\| f\right\| _{L_{p,\lambda } } ,
\]
from here it follows, that
\[
\left\| \tilde{\Im }_{G}^{\alpha } f\right\| _{BMO}
 \le 2\mathop{\sup }\limits_{x,r} \frac{1}{\left|H\left(0,r\right)\right|_{\lambda } } \int\limits _{0}^{r}\left|A_{ch\,\,t}^{\lambda } \tilde{\Im }_{G}^{\alpha } f\left(ch\,\,x\right)-a_{f} \right|sh^{2\lambda }\,\,tdt \lesssim  \left\| f\right\| _{L_{p,\lambda } } .
\]

Theorem 4 is proved.

\textbf{Corollary 3} {\it Let $\alpha p=2\lambda +1,$ $0<\alpha <2\lambda +1,$ $f\in L_{p,\lambda } \left[0,\infty \right).$ If integral $\tilde{\Im }_{G}^{\alpha } f$ is convergence absolutely, then $\tilde{\Im }_{G}^{\alpha } f\in BMO\left[0,\infty \right)$ and the inequality}
\[
\left\| \tilde{\Im }_{G}^{\alpha } f\right\| _{BMO} \lesssim \left\| f\right\| _{L_{p,\lambda}}
\]
is valid.

\

\

\begin{center}
V. S. Guliyev \\
Department of Mathematical Analysis, \\
Institute of Mathematics and Mechanics of NAS of Azerbaijan \\
9, F. Agaev st., Baku, Az 1141, Azerbaijan\\
e-mail: vagif@guliyev.com
\end{center}
\begin{center}
E. J. Ibrahimov\\
Azerbaijan State Oil Academy\\
20, Azadlig Ave., Baku, Az 1010, Azerbaijan\\
e-mail: elmanibrahimov@yahoo.com
\end{center}

\end{document}